\documentclass[a4paper,11pt]{article}
\usepackage{fullpage}
\usepackage{graphics}
\usepackage[latin2]{inputenc}
\usepackage{multicol}
\usepackage{amsmath,amssymb,amsbsy,amsthm,amsfonts}
\usepackage{epsfig} % postscript figures
\begin{document}
%\numberwithin{equation}{section} \marginparwidth=2cm
\def\note#1{\marginpar{\small #1}}

\def\tens#1{\pmb{\mathsf{#1}}}
\def\vec#1{\boldsymbol{#1}}

\def\norm#1{\left|\!\left| #1 \right|\!\right|}
\def\fnorm#1{|\!| #1 |\!|}
\def\abs#1{\left| #1 \right|}
\def\ti{\text{I}}
\def\tii{\text{I\!I}}
\def\tiii{\text{I\!I\!I}}

\def\diver{\mathop{\nabla \cdot }\nolimits}
\def\grad{\mathop{\mathrm{grad}}\nolimits}
\def\Div{\mathop{\mathrm{Div}}\nolimits}
\def\Grad{\mathop{\mathrm{Grad}}\nolimits}

\def\tr{\mathop{\mathrm{tr}}\nolimits}
\def\cof{\mathop{\mathrm{cof}}\nolimits}
\def\det{\mathop{\mathrm{det}}\nolimits}

\def\lin{\mathop{\mathrm{span}}\nolimits}
\def\pr{\noindent \textbf{Proof: }}
\def\pp#1#2{\frac{\partial #1}{\partial #2}}
\def\dd#1#2{\frac{\d #1}{\d #2}}

\def\T{\mathcal{T}}
\def\R{\mathbb{R}}
\def\bx{\vec{x}}
\def\be{\vec{e}}
\def\bef{\vec{f}}
\def\bec{\vec{c}}
\def\bs{\vec{s}}
\def\ba{\vec{a}}
\def\bn{\vec{n}}
\def\bphi{\vec{\varphi}}
\def\btau{\vec{\tau}}
\def\bc{\vec{c}}
\def\bg{\vec{g}}

\def\bW{\vec{W}}
\def\bT{\tens{T}}
\def\bD{\tens{D}}
\def\bF{\tens{F}}
\def\bB{\tens{B}}
\def\bV{\tens{V}}
\def\bS{\tens{S}}
\def\bI{\tens{I}}
\def\bi{\vec{i}}
\def\bv{\vec{v}}
\def\bfi{\vec{\varphi}}
\def\bbsi{\vec{\psi}}
\def\bk{\vec{k}}
\def\b0{\vec{0}}
\def\bom{\vec{\omega}}
\def\bw{\vec{w}}
\def\p{\pi}
\def\bu{\vec{u}}

\def\ID{\mathcal{I}_{\bD}}
\def\IP{\mathcal{I}_{p}}
\def\Pn{(\mathcal{P})}
\def\Pe{(\mathcal{P}^{\eta})}
\def\Pee{(\mathcal{P}^{\varepsilon, \eta})}

\def\Ln#1{L^{#1}_{\bn}}

\def\Wn#1{W^{1,#1}_{\bn}}

\def\Lnd#1{L^{#1}_{\bn, \diver}}

\def\Wnd#1{W^{1,#1}_{\bn, \diver}}

\def\Wndm#1{W^{-1,#1}_{\bn, \diver}}

\def\Wnm#1{W^{-1,#1}_{\bn}}

\def\Lb#1{L^{#1}(\partial \Omega)}

\def\Lnt#1{L^{#1}_{\bn, \btau}}

\def\Wnt#1{W^{1,#1}_{\bn, \btau}}

\def\Lnd#1{L^{#1}_{\bn, \btau, \diver}}

\def\Wntd#1{W^{1,#1}_{\bn, \btau, \diver}}

\def\Wntdm#1{W^{-1,#1}_{\bn,\btau, \diver}}

\def\Wntm#1{W^{-1,#1}_{\bn, \btau}}

%------------------------------------------------
\newtheorem{Theorem}{Theorem}[section]
%{\theorembodyfont{\rmfamily} \newtheorem{Example}{Example}}
\newtheorem{Example}{Example}[section]
\newtheorem{Lem}{Lemma}[section]
\newtheorem{Rem}{Remark}[section]
\newtheorem{Def}{Definition}[section]
\newtheorem{Col}{Corollary}[section]
\newtheorem{Proposition}{Proposition}[section]

%--------------------------------------------------------------------------
\newcommand{\Om}{\Omega}
\newcommand{ \vit}{\hbox{\bf u}}
\newcommand{ \Vit}{\hbox{\bf U}}
\newcommand{ \vitm}{\hbox{\bf w}}
\newcommand{ \ra}{\hbox{\bf r}}
\newcommand{ \vittest }{\hbox{\bf v}}
\newcommand{ \wit}{\hbox{\bf w}}
\newcommand{ \fin}{\hfill $\square$}

\newcommand{\ZZ}{\mathbb{Z}}
\newcommand{\CC}{\mathbb{C}}
\newcommand{\NN}{\mathbb{N}}
\newcommand{\V}{\zeta}
\newcommand{\RR}{\mathbb{R}}
\newcommand{\EE}{\varepsilon}
\newcommand{\Lip}{\textnormal{Lip}}
\newcommand{\XX}{X_{t,|\textnormal{D}|}}
\newcommand{\PP}{\mathfrak{p}}
\newcommand{\VV}{\bar{v}_{\nu}}
\newcommand{\QQ}{\mathbb{Q}}
\newcommand{\HH}{\ell}
\newcommand{\MM}{\mathfrak{m}}
\newcommand{\rr}{\mathcal{R}}
\newcommand{\tore}{\mathbb{T}_3}
\newcommand{\Z}{\mathbb{Z}}
\newcommand{\N}{\mathbb{N}}

\newcommand{\F}{\overline{\boldsymbol{\tau} }}

\newcommand{\moy} {\overline {\vit} }
\newcommand{\moys} {\overline {u} }
\newcommand{\mmoy} {\overline {\wit} }
\newcommand{\g} {\nabla }
\newcommand{\G} {\Gamma }
\newcommand{\x} {{\bf x}}
\newcommand{\E} {\varepsilon}
\newcommand{\BEQ} {\begin{equation} }
\newcommand{\EEQ} {\end{equation} }
\makeatletter
\@addtoreset{equation}{section}
\renewcommand{\theequation}{\arabic{section}.\arabic{equation}}

\newcommand{\hs}{{\rm I} \! {\rm H}_s}
\newcommand{\esp} [1] { {\bf I} \! {\bf H}_{#1} }

\newcommand{\vect}[1] { \overrightarrow #1}

\newcommand{\hsd}{{\rm I} \! {\rm H}_{s+2}}

\newcommand{\HS}{{\bf I} \! {\bf H}_s}
\newcommand{\HSD}{{\bf I} \! {\bf H}_{s+2}}

\newcommand{\hh}{{\rm I} \! {\rm H}}
\newcommand{\lp}{{\rm I} \! {\rm L}_p}
\newcommand{\leb}{{\rm I} \! {\rm L}}
\newcommand{\lprime}{{\rm I} \! {\rm L}_{p'}}
\newcommand{\ldeux}{{\rm I} \! {\rm L}_2}
\newcommand{\lun}{{\rm I} \! {\rm L}_1}
\newcommand{\linf}{{\rm I} \! {\rm L}_\infty}
\newcommand{\expk}{e^{ {\rm i} \, {\bf k} \cdot \x}}
\newcommand{\proj}{{\rm I}Ê\! {\rm P}}

\renewcommand{\theenumi}{\Roman{section}.\arabic{enumi}}

\newcounter{taskcounter}[section]

\newcommand{\bib}[1]{\refstepcounter{taskcounter} {\begin{tabular}{ l p{13,5cm}} \hskip -0,2cm [\Roman{section}.\roman{taskcounter}] & {#1}
\end{tabular}}}

\renewcommand{\thetaskcounter}{\Roman{section}.\roman{taskcounter}}

\newcounter{technique}[section]

\renewcommand{\thetechnique}{\roman{section}.\roman{technique}}

\newcommand{\tech}[1]{\refstepcounter{technique} {({\roman{section}.\roman {technique}}) {\rm  #1}}}

\newcommand{\B}{\mathcal{B}}
%---------------------------------------------------------------------------

\newcommand{\diameter}{\operatorname{diameter}}

%\tableofcontents

%
%
%
%\title[  Critical and subcritical $\alpha$ models of turbulence]{Mathematical results for some $\displaystyle{\alpha}$ models of turbulence with critical and subcritical regularizations }
%
%\author[H. Ali]{Hani Ali}
%\address{IRMAR , UMR CNRS 6625, Universit\'{e} Rennes1, Campus Beaulieu, 35042 Rennes cedex, France}
%\email{hani.ali@univ-rennes1.fr}

%\author[P. Kaplicky]{Petr Kaplicky}
%\address{Charles University, Faculty of Mathematics and
%Physics, Mathematical Institute\\ Sokolovsk\'{a}~83,
%186~75~Prague~8, Czech~Republic}
%\email{Kaplicky@karlin.mff.cuni.cz}

%
%
%\keywords{turbulence model, existence, weak solution}
%\subjclass[2000]{35Q30,35Q35,76F60}

\begin{center}
%\begin{titlepage}
\Large{\textbf{ Global well-posedness  for Deconvolution Magnetohydrodynamics models  with Fractional regularization}}
\end{center}
%\author{\begin{Large}$\hbox{Hani Ali}\thanks{IRMAR, UMR 6625,
%Universit\'e Rennes 1,
%Campus Beaulieu,
%35042 Rennes cedex
%FRANCE;
%hani.ali@univ-rennes1.fr}$ \end{Large}}
%\end{titlepage}
%\date{}
%\date{aujourd'hui}
%\maketitle
\begin{center}
\textbf{{Hani Ali}}$^{\hbox{a}}$\\
$^{\hbox{a}}$\textit{MAP5, CNRS UMR 8145,\\
 Université Paris Descartes,\\
 75006 Paris,
France\\
hani.ali@parisdescartes.fr}
\end{center}
\bigskip
\textbf{Abstract}
\bigskip\\
In this paper, we consider two Approximate Deconvolution Magnetohydrodynamics models which are related to Large Eddy Simulation.
%In this paper, we study  two Approximate Deconvolution Magnetohydrodynamics models which stand for  models of Large Eddy Simulation.
   We first study  existence and uniqueness of solutions  in the double viscous case.  Then, we  study  existence and uniqueness of solutions of the Approximate  Deconvolution MHD model with magnetic diffusivity, but without
kinematic viscosity. In each case, we give the optimal value of regularizations where we can prove
global existence and uniqueness of the solutions.  The second model includes  the Approximate
Deconvolution Euler Model as a particular case. 
%where we study existence and uniqueness of solutions of these two models with optimal regularization value. 
Finally, an asymptotic stability result is shown in the double viscous case with weaker condition on the regularization parameter. 

%  We study   a new regularization of
% the NSE, the   Rotational Approximate Deconvolution model RADM with fractional regularization. We generalize the model studied in Berselli and Lewandowski \cite{} to obtain  the convergence of the solution of the RADM  towards a  solution of the mean Rotational Navier-Stokes equations with weaker  conditions on the parameter regularization.

\smallskip
\smallskip

%\hskip-0.45cm \textbf{Keywords:} {Turbulence simulation and modeling,  Large-eddy simulations, Partial differential equations}\\

%\begin{abstract}
%\end{abstract}
%\maketitle
%\maketitle
%\maketitle
%\begin{multicols}{2}

\section{{Introduction}}
\subsection{Mathematical setting of the problem}
In this paper, we consider the following Approximate Deconvolution MHD equations in a  three dimensional torus $\mathbb{T}_3$

\begin{equation}
\label{dalpha ns}
 \left\{
 \begin{array} {llll} \displaystyle
 \partial_t \vec{w}^{}+ \diver( \overline{ D_{N,\theta}(\vec{w})\otimes   D_{N, \theta}
(\vec{w}^{}) }) -  \diver( \overline{ \vec{B} \otimes  
\vec{B}^{} }) - \nu \Delta \vec{w}^{} + \nabla
q^{} = 0,\\
\displaystyle \partial_t \vec{B}^{}+ \diver( { \vec{B}\otimes   D_{N, \theta}
(\vec{w}^{}) }) -  \diver( { D_{N,\theta}(\vec{w})\otimes 
\vec{B}^{} }) - \mu \Delta \vec{B}^{} = 0,\\
\displaystyle \diver\vec{w}^{} = \diver\vec{B}^{}=0,\
 \displaystyle \int_{\mathbb{T}_3} \vec{w}^{}=  \displaystyle \int_{\mathbb{T}_3} \vec{B}^{}= 0,\\
\vec{w}^{}(t,\vec{x}+{L}\vec{e}_{\vec{j}})=\vec{w}^{}(t,\vec{x}),  \vec{B}^{}(t,\vec{x}+{L}\vec{e}_{\vec{j}})=\vec{B}^{}(t,\vec{x}), \vec{x} \in \mathbb{T}_3,  t >0,\\
\displaystyle \vec{w}^{}_{t=0}=\vec{w}_{0}^{}=\overline{\bv_0}, \vec{B}^{}_{t=0}=\vec{B}_{0}^{},
\end{array}
\right.
\end{equation}
where    the velocity field $\vec{w}^{}(t,\vec{x})$,  the magnetic field $\vec{B}^{}(t,\vec{x})$, and   the scalar pressure $ p(t,\vec{x}) $   are the unknowns, while  $\nu \ge 0 $ is the kinematic viscosity, $\mu \ge 0 $ 
  is the magnetic diffusion, and 
$D_{N,\theta}$  is a deconvolution  operator  given by 
\begin{align}
D_{N, \theta } = \sum_{i = 0}^N ( I - \mathbb{A}_{\theta}^{-1} )^i, \quad \textrm{ for }   0 \le \theta \le 1, \quad \textrm{ and  } N >0,
\end{align}
where
\begin{align}
\label{fractional}
 \mathbb{A}_{\theta}\overline{\bv} = \bv := \left( I + \alpha^{2\theta}(-\Delta)^{\theta} \right) \overline{\bv},
\end{align} 
denotes the Helmholtz
operator with fractional regularization $ \theta$ 
such that  $\bv$ represents the unfiltered velocity and $\overline{\bv}$ is the filtered one. The parameter $\alpha >0$ is the length scale parameter that represents the width of the filter.
In particular, in the filtered equations (\ref{dalpha ns}), the symbol $``\displaystyle{}^{\overline{\quad}}"$  
denotes the above Helmholtz filter (\ref{fractional}), applied component-by-component to the various tensor fields.
If $ \theta=0$,  equations (\ref{dalpha ns}) become the incompressible MHD equations with initial data $(\bv_0, \vec{B}_{0})$. In the absence of the well-posedness theory of the MHD equations, the development of regularized equations is of major importance
for both theoretical and practical purposes. Berselli et al. \cite{bressliduncalewandowski} have suggested two regularizations  for the MHD equations. If we set  $\theta =1$ in (\ref{dalpha ns}), we recover  one of the models  studied  in \cite{bressliduncalewandowski}. For more details about the deconvolution models and the related background, see \cite{AdamStolz, AdamStolzKleiser, dunca06}.\\
The goal of this paper is to establish the existence and the uniqueness of a solution to the system (\ref{dalpha ns}) with fractional regularization $\theta$. First, we will  consider the double viscous  case of (\ref{dalpha ns}) ($ \nu >0$, $ \mu > 0$) in section \ref{sec3} and then  the inviscid case  ($ \nu=0$, $ \mu > 0$) in section \ref{sec4}. In both cases, we will look for the optimal values of $\theta$ where we can
prove global existence and uniqueness of the solutions  for any fixed $N$.

% Before giving our Main results, we first introduce suitable function spaces, then we specify  what we mean by a solution to (\ref{dalpha ns}) in each  cases. 
 Before specifying what we mean by a solution to (\ref{dalpha ns}) in each of these two cases and giving our main results, we introduce suitable
function spaces.

 %Before specifying what we mean by a solution to (\ref{dalpha ns}) in each of these two cases, we introduce suitable
%function spaces.\\
 
\subsection{Functional spaces}
We denote by $L^p(\tore)$ and $H^{s}(\tore)$, ($s \ge -1, \ 1 \le p \le \infty,$) the usual Lebesgue and Sobolev spaces over $\tore$, and define the Bochner spaces $C(0,T;X),$  and  $L^p(0,T;X)$  in the standard way. 
%For the velocity field, we define
The Sobolev spaces $ \vec{H}^{s}=H^s(\tore)^3$,  of mean-free functions are classically characterized in terms of
the Fourier series as follows,
$$\vec{H}^{s} = \left\{ \bv(\bx)=\sum_{{\vec{k} \in \Z^3\setminus\{0\} }}^{}\bc_{\vec{k}} e^{i\vec{k}\cdot \vec{x}}, \left(\bc_{\vec{k}} \right)^{*}=\bc_{-\vec{k}}, \bc_0 = 0,   \|\bv \|_{s,2}^2= \sum_{{\vec{k} \in \Z^3\setminus\{0\} }}^{}| \vec{k} |^{2s} |\bc_{\vec{k}}|^2<\infty  \right\},$$
where  $ \left(\bc_{\vec{k}} \right)^{*}$ denotes the complex conjugate $\bc_{\vec{k}}.$ \\
In addition, we introduce
\begin{align*}
%W^{1,p}_{}&=\left\{ \bv \in W^{s,p}(\tore)^3; \; \int_{\tore}\bv=0  \right\},\\
\vec{H}^{s}_{\sigma }&=\left\{ \bv \in \vec{H}^s; \; \nabla \cdot \bv
=0 \textrm{ in }  \tore \right\},\\
\vec{H}^{-s}&=\left(\vec{H}^{s}_{} \right)^{'},\quad
\vec{L}^2_{}=\vec{H}^{0},\quad 
\vec{L}^2_{\sigma }=\vec{H}^{0}_{\sigma }.
\end{align*}
% Let us mention that by using Poincar\'e inequality we have 
%  \begin{align}
% \|\vec{v} \|_{s,2} \approx \| \widetilde{\vec{v}}\|_{s+2\theta_1,2} \approx \| \overline{\vec{v}} \|_{s+2\theta_2,2}. 
% \end{align}
%For $s  \in \R$,
%let $\vec{I}=\left\lbrace \bk  \ \hbox{such  that}  \   \displaystyle \bk=\frac{2\pi \ba}{\vec{L}},\  \ba \in \Z^{3} , \ba \neq 0 \right\rbrace ,$\\
%the usual Sobolev spaces $W^{s,2}$  with zero space average are classicaly characterized in terms of Fourier series  
%$$W^{s,2}=\{\bv=\sum_{\bk \in \vec{I}}\hat{\bv}_{\bk}\exp{\left\lbrace   i\bk\cdot \bx\right\rbrace }, \
%\| \bv \|_{{s,2}}^{2}< \infty          \},$$
%where 
%$$\| \bv \|_{{s,2}}^{2}= \sum_{\bk \in \vec{I}}|\bk|^{2s} |\hat{\bv}_{\bk}|^{2}.$$
%Therefore,  the space $W^{s,2}_{\diver}$ is  
%$$W^{s,2}_{\diver}=\{\bv=\sum_{\bk \in \vec{I}}\hat{\bv}_{\bk}\exp{\left\lbrace   i\bk\cdot \bx\right\rbrace }, \ 
%\ \vec{k}\cdot\hat{\bv}_{\bk}=0, \
%\| \bv \|_{{s,2}}^{2}< \infty          \}.$$

\subsection{Main results}
In this subsection, we define the notion of a solution to the system (\ref{dalpha ns}) and we give precisely the main results of this paper.
\subsubsection{Definitions of regular weak solutions}
 
We first consider the double viscous  case of (\ref{dalpha ns}) ($ \nu >0$, $ \mu > 0$).

\begin{Def} 
\label{definition1}
Let  $ \nu >0$, $ \mu > 0$, $\bv_0 \in L^{2}_{\sigma }$ and $\vec{B}_0 \in L^{2}_{\sigma}$.  For any $0 \le \theta \le 1$ 	and $ 0 \le N <\infty  $ we say that the triplet $(\bw, \vec{B}, q) $  is a `` regular" weak solution to (\ref{dalpha ns}) if for any $T>0$ the following
properties are satisfied:
\begin{align}
\bw &\in \mathcal{C}_{}(0,T;\vec{H}^{\theta}_{\sigma }) \cap
L^2(0,T;\vec{H}^{1+\theta}_{\sigma }), \quad
\vec{B} \in \mathcal{C}_{}(0,T;\vec{L}^{2}_{\sigma }) \cap
L^2(0,T;\vec{H}^{1}_{\sigma }), \label{bv12}\\
\partial_t\bw&\in   L^{2}(0,T;\vec{H}^{2 \theta - \frac{3}{2}}_{}), \quad 
\partial_t\vec{B} \in   L^{2}(0,T;\vec{H}^{-1}_{}),
 \label{bvt} \\
q&\in L^{2}(0,T;{H}^{2 \theta - \frac{1}{2}}(\tore))
\label{psp},
\end{align}
the triplet $(\bw, \vec{B}, q)$ fulfill
\begin{equation}
\begin{split}
\int_0^T \langle \partial_t\bw, \bfi \rangle  -  \langle  \overline{D_{N,\theta}(\bw) \otimes  D_{N,\theta}(\bw)}, 
\nabla \bfi \rangle  +  \nu \langle 
\nabla \bw, \nabla \bfi  \rangle  +\langle \nabla q,  \bfi \rangle  \; dt\\
  + \int_0^T  \langle  \overline{(\vec{B}) \otimes (\vec{B})}, \nabla \bfi \rangle \; dt = 0 
\qquad \textrm{ for all } {\bfi}\in L^{2}(0,T; \vec{H}^{ \frac{3}{2} -2 \theta }),
 %\vec{H}^{\frac{3}{2}-2\theta}_{}\cap  \vec{H}^{1-\theta}_{}\cap \vec{H}^{1-2\theta}_{}
\end{split}\label{weak1}
\end{equation}
\begin{equation}
\begin{split}
\int_0^T \langle \partial_t\vec{B}, \bfi \rangle  +  \langle  D_{N,\theta}(\bw) \otimes  \vec{B}, 
\nabla \bfi \rangle  +  \mu  \langle 
\nabla \vec{B}, \nabla \bfi  \rangle  \; dt \qquad  \qquad \qquad \\
    - \int_0^T  \langle \vec{B} \otimes D_{N,\theta}(\bw), \nabla \bfi \rangle \; dt = 0
\qquad \textrm{ for all } {\bbsi}\in L^{2}(0,T; \vec{H}^{1}).
 %\vec{H}^{\frac{3}{2}-2\theta}_{}\cap  \vec{H}^{1-\theta}_{}\cap \vec{H}^{1-2\theta}_{}
\end{split}\label{weakb1}
\end{equation}
Moreover,
\begin{equation}
\bw(0)= \bw_0 \quad \textrm{ and } \quad \vec{B}(0)= \vec{B}_0.
\label{intcond}
\end{equation}
\end{Def}
 Then, we give the definition of the solution of  (\ref{dalpha ns}) with partial viscous term ($ \nu =0$, $ \mu > 0$).

\begin{Def} 
Let  $ \nu =0$, $ \mu > 0$, $\bv_0 \in \vec{L}^{2}_{\sigma}$ and $\vec{B}_0 \in \vec{L}^{2}_{\sigma}$.  For any $0 \le \theta \le 1$ 	and $ 0 \le N <\infty  $ we say that the triplet $(\bw, \vec{B}, q) $  is a  weak solution to (\ref{dalpha ns}) if for any $T>0$ the following
properties are satisfied:
\begin{align}
\bw &\in L^{\infty}(0,T;\vec{H}^{\theta}_{\sigma }) \cap \mathcal{C}_{}(0,T;\vec{L}^{2}_{\sigma }),\label{nubv12}\\
\vec{B} &\in \mathcal{C}_{}(0,T;\vec{L}^{2}_{\sigma }) \cap
L^2(0,T;\vec{H}^{1}_{\sigma }),\label{nuB12}\\
\partial_t\bw&\in   L^{4}(0,T;\vec{H}^{2 \theta - {2}}_{})\cap L^{\frac{6}{5}}(0,T ; \vec{H}^{2 \theta - \frac{5}{6}}),
\label{nubvt}\\
\partial_t\vec{B}&\in   L^{2}(0,T;\vec{H}^{-1}_{}),
\label{nuBt}\\
q&\in L^{2}(0,T;{H}^{2 \theta -{1}}(\tore))
\label{nupsp},
\end{align}
the triplet $(\bw, \vec{B}, q)$ fulfill
\begin{equation}
\begin{split}
\int_0^T \langle \partial_t\bw, \bfi \rangle  -  \langle  \overline{D_{N,\theta}(\bw) \otimes  D_{N,\theta}(\bw)}, 
\nabla \bfi \rangle  +\langle \nabla q,  \bfi \rangle  \; dt\\
  + \int_0^T  \langle  \overline{(\vec{B}) \otimes (\vec{B})}, \nabla \bfi \rangle \; dt = 0 
\qquad \textrm{ for all } {\bfi}\in L^{\frac{4}{3}}(0,T; \vec{H}^{2 - 2 \theta}),
 %\vec{H}^{\frac{3}{2}-2\theta}_{}\cap  \vec{H}^{1-\theta}_{}\cap \vec{H}^{1-2\theta}_{}
\end{split}\label{nuweak1}
\end{equation}
\begin{equation}
\begin{split}
\int_0^T \langle \partial_t\vec{B}, \bfi \rangle  +  \langle  D_{N,\theta}(\bw) \otimes  \vec{B}, 
\nabla \bfi \rangle  +  \mu  \langle 
\nabla \vec{B}, \nabla \bfi  \rangle  \; dt \qquad  \qquad \qquad \\
    - \int_0^T  \langle \vec{B} \otimes D_{N,\theta}(\bw), \nabla \bfi \rangle \; dt = 0
\qquad \textrm{ for all } {\bbsi}\in L^{2}(0,T; \vec{H}^{1}).
 %\vec{H}^{\frac{3}{2}-2\theta}_{}\cap  \vec{H}^{1-\theta}_{}\cap \vec{H}^{1-2\theta}_{}
\end{split}\label{nuweakb1}
\end{equation}
Moreover,
\begin{equation}
\bw(0)= \bw_0 \quad \textrm{ and } \quad \vec{B}(0)= \vec{B}_0.
\label{intialecondition}
\end{equation}
Furthermore,
if,  $ \vec{B}_0 \in \vec{H}^{1}_{\sigma}$  we say that  $(\bw, \vec{B}, q) $  is a ``regular"  weak  solution to (\ref{dalpha ns}) if it is a weak solution, and in addition $ \bw$ and  $\vec{B}$ satisfy
\begin{align}
\bw &\in \mathcal{C}_{}(0,T;\vec{H}^{\theta}_{\sigma }) \cap L^{\infty}(0,T;\vec{H}^{2\theta}_{\sigma })  ,\label{bv12bv12}\\
\vec{B} &\in \mathcal{C}_{}(0,T;\vec{H}^{1}_{\sigma }) \cap
 L^2(0,T;\vec{H}^{2}_{\sigma }).\label{B12}
%\\
% \bw_{,t}&\in   L^{2}(0,T;\vec{H}^{4 \theta - \frac{5}{2}}_{}),
% \label{bvt}\\
% \vec{B}_{,t}&\in   L^{2}(0,T;\vec{H}^{-1}_{}),
% \label{Bt}
\end{align}
 
\end{Def}
Note that  the term ``regular"  in the above definitions  is used when weak solutions are unique and do not develop a finite-time singularities.

\subsubsection{Statement of the main results}
The main results of this paper are the following:

% The first main purpose of this paper is to investigate the double viscous case of (\ref{dalpha ns}), We will prove the following.
 \begin{Theorem} 
Consider the Approximate Deconvolution MHD equations (\ref{dalpha ns})  with $ \nu >0$ and  $ \mu > 0.$  Assume   $\bv_0 \in L^{2}_{\sigma }$ and $\vec{B}_0 \in \vec{L}^{2}_{\sigma }$.  Let $ \theta \ge \frac{1}{2}$ 	and  let $ 0 \le N <\infty  $ is given and fixed.  Then  (\ref{dalpha ns})  with initial data $(\bv_0,\vec{B}_{0})$  has a unique regular
weak solution.
\label{TH1}
\end{Theorem}
% The second main purpose of this paper is to investigate the above Deconvolution MHD system with partial viscous term $ \nu =0$ and  $ \mu > 0.$  We will prove the following.
 \begin{Theorem}
 \label{Theorem2}
Consider the Approximate Deconvolution MHD equations (\ref{dalpha ns})  with $ \nu =0$ and  $ \mu > 0.$  Assume   $\bv_0 \in L^{2}_{\sigma }$ and $\vec{B}_0 \in L^{2}_{\sigma }$.  Let $ \theta \ge \frac{5}{6}$ 	and  let $ 0 \le N <\infty  $ is given and fixed.  Then  (\ref{dalpha ns})  with initial data $(\bv_0,\vec{B}_{0})$ has a 
weak solution.
Furthermore, if $\vec{B}_0 \in \vec{H}^{1}_{\sigma }$  then   (\ref{dalpha ns})  with initial data $(\bv_0,\vec{B}_{0})$ has a unique ``regular''  weak solution among the class of weak solutions.
\end{Theorem}

Theorem \ref{Theorem2} includes a special case that deserves a separate formulation.
If we consider the case $ \vec{B} =0$, then equations (\ref{dalpha ns}) with $ \nu =0$
  become    the Approximate Deconvolution Model of  Euler equations  (ADE)  with periodic 
boundary conditions. Consequently, Theorem \ref{Theorem2}  reduces to the following statement.

\begin{Col}
 \label{Col2}
  Assume   $\bv_0 \in L^{2}_{\sigma }$.  Let $ \theta \ge \frac{5}{6}$ 	and  $ 0 \le N <\infty  $ given and fixed numbers.  Then  the Approximate Deconvolution Euler equations  (ADE)   with initial data $\bv_0$ have a unique regular weak solution $
\bw \in \mathcal{C}_{}(0,T;\vec{H}^{\theta}_{\sigma }).$  This regular weak solution satisfies the energy equality:
  \begin{equation}
\begin{array}{lll}
 \label{energyequality}
\displaystyle \frac{1}{2}\| \mathbb{A}_{\theta}^{\frac{1}{2}}D_{N,\theta}^{\frac{1}{2}}({\bw})(t) \|_{2}^2 
= \displaystyle \frac{1}{2}\| \mathbb{A}_{\theta}^{\frac{1}{2}}D_{N,\theta}^{\frac{1}{2}}(\overline{\bv}_{0})\|_{2}^2.
\end{array}
\end{equation}
\end{Col}

\begin{Rem}
From the above energy equality, we can deduce by using the properties of $D_{N,\theta}$ and $\mathbb{A}_{\theta}$ (see section \ref{decon}), the following inequality,
 \begin{equation}
\begin{array}{lll}
 \label{energyequalityalpha}
\displaystyle \| {\bw}(t) \|_{2}^2 +  \alpha^{2\theta}\| {\bw}(t) \|_{\theta,2}^2
 \le \displaystyle \| {\bv}_{0} \|_{2}^2.
\end{array}
\end{equation}
Using the same techniques as in  \cite{Lariostiti}, we can derive from  (\ref{energyequalityalpha}) the following criterion
for finite-time blow-up of the Euler equations.\\ 
Blow-up criterion: Assume ${\bv}_{0} \in \vec{H}_{\sigma}^s$, for some $s > 3$. Suppose that
there exists a finite time $T^* > 0$ such that the solution ${\bw}$  of (ADE) with  initial data ${\bv}_{0}$    satisfies the below inequality  for each $\alpha > 0$  	and $ 0 \le N <\infty$: 
 \begin{equation}
\begin{array}{lll}
 \label{blowup}
\displaystyle \sup_{t \in [0,T^*]} \limsup_{\alpha \rightarrow 0}  \alpha^{2\theta}\| {\bw}(t) \|_{\theta,2}^2 > 0.
\end{array}
\end{equation}
Then the Euler equations with initial data ${\bv}_{0}$  develop a singularity in the interval
$ [0, T^*].$
\end{Rem}

Let us  briefly review some existence results for  regularizations of the magnetohydrodynamics equations; we will not attempt to address exhaustive references in this paper.\\
% When $N=0$ (\ref{dalpha ns}) is similar to the MHD-voight. The voight regularization is studied first in the context of the Navier-Stokes equations in  \cite{}. 
 The existence and  uniqueness of the inviscid Voigt regularization of the MHD equations are studied in \cite{LariostitiMHD,CaSe11}.
When $\theta=1 $ and $N$ fixed, the global well-posedness of the double viscous problem has been proven in \cite{bressliduncalewandowski}. When $ \vec{B} = 0$ and $\theta=1 $, equations (\ref{dalpha ns})  become the well-known Approximate Deconvolution Model (ADM).  Recently, the authors
\cite{bresslilewandowski} have proven that the ADM with fractional regularization has a unique regular weak solution when $\theta > \frac{3}{4}$. 
The question of the limit behavior of the ADM solutions when $N$ tends to infinity is studied also in \cite{bressliduncalewandowski} when $\theta > \frac{3}{4}$. The above uniqueness and   convergence  results are improved in \cite{HaniAli} to cover the range when  $ \theta \ge \frac{1}{6}$. 
 The deconvolution operator, for different values of $\theta$,  was used in \cite{DuncaLewandowski} to study  the  rate of convergence of the ADM model to the mean Navier-Stokes Equations. 
Very recently, Berselli et al.  \cite{bressliduncalewandowski} have adapted the results  in \cite{bresslilewandowski} to the  Deconvolution MHD equations (\ref{dalpha ns}).

In this paper, our approach is partially inspired from the recent studies  of Berselli et al.  \cite{bressliduncalewandowski,bresslilewandowski}.  We point out that in \cite{bressliduncalewandowski} the authors consider only the double viscous case with $\theta=1$. In this paper, we study  not only the double viscous case but also the partial inviscid case with fractional regularization. In the double viscous case, we derive some apriori estimates that are uniform with respect to $N$. These uniform estimates combined with a compactness argument allow us to study  the limit behavior of the solutions of the Deconvolution MHD model when $N$ tends to infinity, for $ \theta \ge \frac{1}{2}$. Unfortunately, in the partial viscous case,  we  are not able to use the method  developed in this paper
to study the limit behavior of the solutions when $N$ tends to infinity. This question remains open. \\

The remaining part of the paper is organized as follows. In section 2, we start by giving some preliminaries
about  the deconvolution operator $D_{N,\theta}$.  Section 3 provides  existence and uniqueness results in the double viscous case  where $  0 \le N < \infty$ is fixed and $\theta \ge \frac{1}{2}$. Section 4 provides  existence and uniqueness results in the partial inviscid  case  where $  0 \le N < \infty$ is fixed  and $\theta \ge \frac{5}{6}$.    The last section provides a result  of consistency and stability   which illustrates   how the system approximates the  MHD equations.

\section{The deconvolution operator}
\label{decon}

In this paper we consider a  generalized deconvolution operator  \cite{bresslilewandowski,A01}. 
%The deconvolution operator considered in this paper is a generalized deconvolution operator  \cite{bresslilewandowski,A01}. 
This deconvolution operator is  constructed by using $\displaystyle \mathbb{A}_{\theta}$, the Helmholtz equation with fractional regularization.  

For $ p = -1, -\frac{1}{2}, \frac{1}{2} \textrm { or } 1  $, the non-local operator   $\mathbb{A}_{\theta}^{p}$    is  defined in the periodic case as
\BEQ 
\label{2.1} 
\mathbb{A}_{\theta}^{p}  \left ( \sum_{ {\bk}\in \Z^3\setminus\{0\}} \bc_{\bk} \expk \right ) =  
\sum_{ {\bk}\in \Z^3\setminus\{0\}}\left(1+ \alpha^{2 \theta} | {\bk} |^{2 \theta}\right)^{p}   \bc_{\bk} \expk, \textrm{ for any  }  0 \le \theta \le 1, \EEQ 
 and in particular  $\mathbb{A}_{\theta}^{-1}$  verifies 
 \BEQ
\label{lemme2.2}
\| \mathbb{A}_{\theta}^{-1}({\bv})\|_{s+2\theta,2}  \le {\alpha^{-2\theta}} \|\bv \|_{s,2}.
\EEQ
  This operator is a
differential filter \cite{germano}, that commutes with differentiation under periodic boundary conditions.
In the analysis we will need the following relations  which follow from (\ref{2.1})
 \BEQ
\label{lemmeatheta1}
\left(\mathbb{A}_{\theta}^{-\frac{1}{2}}(\bv), \mathbb{A}_{\theta}^{\frac{1}{2}}(\bw)\right) = \left(\bv,\bw\right), \textrm{ for any } \bv \in \vec{L}^{2} \textrm{ and } \bw \in \vec{H}^{\theta},
\EEQ
\BEQ
\label{lemmeatheta2}
 \| \mathbb{A}_{\theta}(\bv)\|_{2 }^2=\| \bv\|_{2 }^2 + 2\alpha^{2\theta} \| \bv\|_{\theta,2 }^2 + \alpha^{4\theta}\| \bv\|_{2\theta,2 }^2, \textrm{ for any } \bv \in   \vec{H}^{2\theta}.
\EEQ
%{\bf Proof}. By definition, one has 
%$$  || \moys - u ||_{s,2}^2 = \sum_{ {\bk}Ê\in {\cal  I}_3} \left ( {\alpha^2 |  {\bk}Ê |^2 
%\over 1 + \alpha^2 |  {\bk}Ê |^2 } \right ) ^2 |  {\bk}Ê |^{2 s} | u_{\bk}Ê|^2.   $$
%Let $\E >0$. As $u \in \hs$, there exists $N$ be such that 
%$$ \sum_{{\bk} \in  {\cal  I}_3 \setminus I_N} | {\bk}Ê|^{2s} | u_{\bk} |^2 \le {\E \over 2}, $$
%and since $ {\alpha^2 |  {\bk}Ê |^2 
%/ (1 + \alpha^2 |  {\bk}Ê |^2) } \le 1$, 
%$$ {\bf I} _N = \sum_{ {\bk} \in  {\cal  I}_3 \setminus I_N} \left ( {\alpha^2 |  {\bk}Ê |^2 
%\over 1 + \alpha^2 |  {\bk}Ê |^2 } \right ) ^2 |  {\bk}Ê |^{2 s} | u_{\bk}Ê|^2 \le {\E \over 2}.$$
%On the other hand, because the set $I_N$ is finite, 
%$$ \lim_{ \alpha \rightarrow 0} \sum_{ {\bk} \in   I_N^\star} \left ( {\alpha^2 |  {\bk}Ê |^2 
%\over 1 + \alpha^2 |  {\bk}Ê |^2 } \right ) ^2 |  {\bk}Ê |^{2 s} | u_{\bk}Ê|^2 = 0.$$
%Therefore, there exists $\alpha_0 >0$ be such that for each $\alpha \in \, ]0, \alpha_0[$ 
%one has 
%$${\bf J}_N = \sum_{ {\bk} \in   I_N^\star} \left ( {\alpha^2 |  {\bk}Ê |^2 
%\over 1 + \alpha^2 |  {\bk}Ê |^2 } \right ) ^2 |  {\bk}Ê |^{2 s} | u_{\bk}Ê|^2 \le {\E \over 2}.$$
%As $ || \moys - u ||_{s,2}^2  = {\bf I}_N + {\bf J}_N$, then for all $\alpha \in \, ]0, \alpha_0[$, one has $ || \moys %- u ||_{s,2}^2  \le \E$ ending the proof like that. 

A straightforward calculation yields

\BEQ D_{N,\theta}  \left ( \sum_{ {\bk}\in \Z^3\setminus\{0\}} \bc_{\bk} \expk \right ) =  
\sum_{ {\bk}\in \Z^3\setminus\{0\}}\left(1+ \alpha^{2 \theta} | {\bk} |^{2 \theta}\right) \left ( 1- \left ( \frac{\alpha^{2 \theta} | {\bk} |^{2 \theta}}
 {1 +  \alpha^{2 \theta} | {\bk} |^{2 \theta}}\right )^{N+1} \right )  \bc_{\bk} \expk. \EEQ 
Thus 
 \BEQ D_{N,\theta}  \left ( \sum_{ {\bk}\in \Z^3\setminus\{0\}} \bc_{\bk} \expk \right ) =  \sum_{ {\bk}\in \Z^3\setminus\{0\}}  \widehat{D_{N,\theta}}(\bk)  \bc_{\bk} \expk, \EEQ 
 where we have for all ${\bk}\in {\cal  I}_3,$
 
 \begin{align}
 \widehat{D_{0,\theta}}(\bk) &= 1, \\
 1 \le  \widehat{D_{N,\theta}}(\bk) &\le N+1 \quad \hbox { for each } N >0, \\
 \hbox{ and  }   \widehat{D_{N,\theta}}(\bk) &\le \widehat{\mathbb{A}_{\theta}} := \left(1+ \alpha^{2 \theta} | {\bk} |^{2 \theta}\right) \hbox{ for  a fixed } \alpha >0. \label{haniali}
  \end{align}

The following elementary lemma is given in \cite{bresslilewandowski} and  will play an important role.

\begin{Lem}
\label{lemme1}
For all  $s \ge -1$, ${\bk}\in \Z^3\setminus\{0\},$ and for each $N >0$ there exist a constant $C >0$ such that for all $\bv$ sufficiently smooth  we have 
\begin{align}
  &\|\bv\|_{s,2} \le    \|D_{N,\theta} \left ( \bv \right )\|_{s,2} \le  (N+1)  \|\bv\|_{s,2},\label{lemme13a}\\
    &\|\bv\|_{s,2} \le   C \|D_{N,\theta} \left ( \bv \right )\|_{s,2} \le  C \|{\mathbb{A}_{\theta}}^{\frac{1}{2}}D_{N,\theta}^{\frac{1}{2}}(\bv)\|_{s,2},\label{lemme13b}\\
     &\|{\mathbb{A}_{\theta}}^{\frac{1}{2}}D_{N,\theta}^{\frac{1}{2}}(\overline{\bv})\|_{s,2}  \le  \|\bv\|_{s,2}, \label{lemme13}\\    
      &\|\bv\|_{s+\theta,2} \le   C(\alpha) \|{\mathbb{A}_{\theta}}^{\frac{1}{2}}D_{N,\theta}^{\frac{1}{2}}(\bv)\|_{s,2}.\label{lemme13d}
  \end{align}
\end{Lem}

% Therefore let $ (\mathbb{A}_{\theta} \bw,\mathbb{A}_{\theta} q)$  is a  distributional solution to the rotational  Navier-Stokes equations  (\ref{nsBM})-(\ref{nsBLM})

%\section{The double viscous case}

\section{Regular weak solution in the double viscous case} 
\label{sec3}

In this section we prove Theorem \ref{TH1} by using the Galerkin method.
%The proof of Theorem \ref{TH1} follows the classical scheme. 
For simplicity  we will consider only the case when $0<N<\infty $ and the case $N=0$ is deduced  by replacing $D_{0,\theta}$ by $I$. 
 We start by constructing
approximated solutions $(\bw^n, \vec{B}^n, q^n)$ via Galerkin method. Then we seek for a priori
estimates that are uniform with respect to $n$. Next, we passe to the limit in the
equations after having used compactness properties. Finally we show that the solution
we constructed is unique thanks to Gronwall's lemma \cite{RT83}.\\

\textbf{Step 1}(Galerkin approximation).
Consider a sequence $\left\{ \bfi^{r} \right\}_{r=1}^{\infty}$ consisting of $L^2$-orthonormal  and $H^{1}$-orthogonal eigenvectors of the Stokes problem subjected to the
space periodic conditions. 
%\begin{align}
%-\Delta \bfi^{r} = |r|^2 \bfi^{r},\hbox{ in } \tore \quad \int_{\tore}\bfi^{r}=0,  \quad {\bfi^{r}}( \bx + L \vec{e}_{\vec{j}})={\bfi^{r}}(\bx) \hbox{ for all } {r} \in \Z^3\setminus\{0\} .
%\end{align}
We note that this sequence forms a hilbertian basis of $L^2$.\\
%Stokes operator subjected to the space periodic conditions (\ref{bc1}) and 
We set 
\begin{equation}
\begin{split}
\bw^n(t,\bx)=\sum_{{r=1}}^{n}\bc_{r}^n (t) \bfi^{r}(\bx),\\
\vec{B}^n(t,\bx)=\sum_{{r=1}}^{n}\vec{d}_{r}^n (t) \bfi^{r}(\bx),
\\
\quad  \hbox{ and } q^n(t,\bx)=\sum_{{|\vec{k}|=1}}^{n}q_{\vec{k}}^n (t) e^{i \vec{k} \cdot \bx}.
\end{split}
\end{equation}
% such that $\bk \cdot  \bc_{r}^n = 0  $ for all $\bk \in \Z^3\setminus\{0\}$ and $\left(\bc_{r}^n \right)^{*}=\bc_{-r}^n  $, where  $ \left(\bc_{r}^n \right)^{*}$ denote the complex conjugate $\bc_{r}^n.  $
% Thus due of (\ref{TKEbis}) and (\ref{TKE}) we have 
%\begin{equation}
% \widetilde{\bw}^n(t,\bx)=\sum_{{|r|=1}}^{N}\widetilde{\bc}_{r}^n (t) \bfi^{r}(\bx) \textrm{ and }  \overline{\bw}^n(t,\bx)=\sum_{{|r|=1}}^{N}\overline{\bc}_{r}^n (t) \bfi^{r}(\bx), 
%\end{equation}
% where 
%\begin{equation}
% \widetilde{\bc}_{r}^n = \frac{{\bc}_{r}^n }{1+\alpha^{2\theta_1}|r|^{2\theta_1}} \textrm{ and }  \overline{\bc}_{r}^n = \frac{{\bc}_{r}^n }{1+\alpha^{2\theta_2}|r|^{2\theta_2}}, 
%\end{equation}
%for all $\bk \in \Z^3\setminus\{0\}$.\\
We look for $(\bw^n(t,\bx),\vec{B}^n(t,\bx), q^n(t,\bx)) $ that are determined through the system of equations 
\begin{equation}
\begin{split}
 \left( \partial\bw^n, \bfi^{r} \right)  -  \left(\overline{D_{N,\theta}(\bw^n) \otimes  D_{N,\theta}(\bw^n)} , 
\nabla\bfi^{r}\right) +  \nu \left(
\nabla \bw^n, \nabla \bfi^{r}  \right)\; \\
+ \left( \overline{\vec{B}^n \otimes  \vec{B}^n} , \nabla \bfi^{r} \right) = 0 \; ,
\quad {r=1},2,...,n,
\end{split}\label{weak1galerkine}
\end{equation}
\begin{equation}
\begin{split}
 \left( \partial_t\vec{B}^n, \bfi^{r} \right)  +  \left({D_{N,\theta}(\bw^n) \otimes  \vec{B}^n} , 
\nabla\bfi^{r}\right) +  \mu \left(
\nabla \vec{B}^n, \nabla \bfi^{r} \right)\; \\
- \left( {\vec{B}^n \otimes  D_{N,\theta}(\bw^n)} , \nabla \bfi^{r} \right) = 0 \; ,
\quad {r=1},2,...,n,
\end{split}\label{weakb1galerkine}
\end{equation}

and 
%\begin{equation}
%\begin{split}
% q^n = - \sum_{i,j}\partial_i\partial_j \Delta^{-1}(\Pi^n(\overline{{v}_i^n {v}_j^n)})= - \sum_{i,j} R_{ij}(\Pi^n(\widetilde{v}_i^n\overline{v}_j^n)),
%\end{split}\label{pressuregalerkine}
%\end{equation}
\begin{equation}
\begin{split}
 %\displaystyle \Delta^{} q^n=  \nabla \cdot   \Pi^n \left(  \overline{D_{N} {\bw}^n  \times \nabla \times D_{N} {\bw}^n }  \right).
\displaystyle \Delta^{} q^n=  -\diver \diver  \left(\Pi^n(\overline{D_{N,\theta}({\bw}^n) \otimes D_{N,\theta}({\bw}^n)} - \overline{\vec{B}^n \otimes  \vec{B}^n})  \right).
\end{split}\label{pressuregalerkine}
\end{equation}
Where the projector $ \displaystyle \Pi^n $ assign to any Fourier series $\displaystyle \sum_{\bk \in \Z^3\setminus\{0\}} \vec{g}_{\bk} e^{i\bk \cdot\bx} $ the following  series  
 $\displaystyle \sum_{\bk \in \Z^3\setminus\{0\}, |\bk| \le n} \vec{g}_{\bk} e^{i\bk \cdot\bx}. $
 %we deduce from the elliptic equation (\ref{pressuregalerkine}) that 
 % after dropping indices of $N$ for shortness, that 
  % \begin{equation}
%\begin{split}
% p^n = - \sum_{i,j}\partial_i\partial_j \Delta^{-1}(\Pi^n(\widetilde{v}_i^n\overline{v}_j^n))= - \sum_{i,j} R_{ij}(\Pi^n(\widetilde{v}_i^n\overline{v}_j^n)),
%\end{split}\label{pressurepseudo}
%\end{equation}
% and $R_{ij}$ is the Riez operator defined through the Fourier transform by
% \begin{equation}
%\begin{split}
%  \widehat{R_{ij}(u)}=\frac{k_i k_j}{|\bk|^2}\widehat{u({\bk})}, \quad \hbox{ for all } \bk \in \Z^3\setminus\{0\}.
%  \end{split}\label{pressurepseudo2}
%\end{equation}

% The Laplace operator $\displaystyle \Delta^{-1}$ can be  defined through the Fourier transform hence $\displaystyle (\Delta^{-1}\diver \diver)$ may be viewed as a pseudofifferentiel operator 
%  \begin{equation}
%\begin{split}
%   ((\Delta^{-1}\diver \diver){\bw})_{\bk}=\frac{\bk^2}{|\bk|^2}\bw_{\bk}.
%\end{split}\label{pressurepseudo}
%\end{equation}
Moreover we require that $\bw^n $ and $\vec{B}^n $  satisfy the following initial conditions
\begin{equation}
\label{initial Galerkine}
\bw^n(0,.)= \bw^n_0= \sum_{r=1}^{n}\bc_0^n  \bfi^{r}(\bx), \quad  \vec{B}^n(0,.)= \vec{B}^n_0= \sum_{r=1}^{n}\vec{d}_0^n  \bfi^{r}(\bx)
\end{equation}
and 
\begin{equation}
\begin{split}
 \label{initial2 Galerkine}
\bw^n_0 \rightarrow \bw_0  \quad \textrm{ strongly  in }  \vec{H}^{\theta}_{\sigma} \quad \textrm{ when } n \rightarrow \infty,\\
\vec{B}^n_0 \rightarrow \vec{B}_0  \quad \textrm{ strongly  in }  \vec{L}^{2}_{\sigma} \quad \textrm{ when } n \rightarrow \infty.
\end{split}
\end{equation}
%Where the initial condition $ \overline{\bw}^n_0$ is deduced from $\bw^n_0$ through the relation (\ref{TKE}).
%\begin{equation}
% \label{initialbar Galerkine}
%  \overline{\bw}^n_0   + \alpha^{2\theta}(-\Delta)^{\theta} \overline{\bw}^n_0=\bw^n_0
%\end{equation}
 The classical Caratheodory theory \cite{Wa70}  then implies the short-time existence of solutions 
to (\ref{weak1galerkine})-(\ref{pressuregalerkine}).  Next we derive  estimates on $\bc^n$ and $\vec{d}^n$ that are uniform w.r.t. $n$.
These estimates then imply that the  solution of  (\ref{weak1galerkine})-(\ref{pressuregalerkine}) constructed on a short time interval $[0, T^n[ $ exists for all $t \in [0, T]$.\\

\textbf{Step 2} (A priori estimates)
Multiplying the $r$th equation in (\ref{weak1galerkine}) with $\alpha^{2\theta}|\bk|^{2\theta}\widehat{D_{N,\theta}}{\bc}^n_{r}(t)+\widehat{D_{N,\theta}}{\bc}^n_{r}(t)$,  and the $r$th equation in (\ref{weakb1galerkine}) with ${\vec{d}}^n_{r}(t)$ summing over ${r=1},2,...,n$, integrating over time from $0$ to $t$ and using the following identities
 \begin{equation}
 \label{divergencfreebar1}
 \begin{split}
\left(\partial_t{\bw}^n, \mathbb{A}_{\theta}D_{N,\theta}({\bw}^n)\right)=
\frac{1}{2}\frac{d}{dt}\|\mathbb{A}_{\theta}^{\frac{1}{2}}D_{N,\theta}^{\frac{1}{2}}({\bw}^n) \|_{2}^2, \\
\left(\partial_t{\vec{B}}^n, {\vec{B}}^n\right)=
\frac{1}{2}\frac{d}{dt}\|{\vec{B}}^n \|_{2}^2,
\end{split}
\end{equation}
\begin{equation}
 \label{divergencfreebar2divergencfreebar2}
 \begin{split}
\left(-\Delta{\bw}^n,   \mathbb{A}_{\theta}^{}D_{N,\theta}^{}({\bw}^n)\right)=
\|\mathbb{A}_{\theta}^{\frac{1}{2}}D_{N,\theta}^{\frac{1}{2}}({\bw}^n)\|_{1,2}^2,\\
\left(-\Delta{\vec{B}}^n,   {\vec{B}}^n\right)=
\|{\vec{B}}^n\|_{1,2}^2,
\end{split}
\end{equation}
% \begin{equation}
%  \label{divergencfreebar3}
% \langle \overline{\bef}, \mathbb{A}_{\theta}^{}D_{N,\theta}^{}{\bw}^n \rangle= \langle\mathbb{A}_{\theta}^{\frac{1}{2}}D_{N,\theta}^{\frac{1}{2}} \overline{\bef}, \mathbb{A}_{\theta}^{\frac{1}{2}}D_{N,\theta}^{\frac{1}{2}}{\bw}^n \rangle,
% \end{equation}

\begin{equation}
 \label{divergencfreebar}
 \begin{array}{lll}
\left(\overline{{D_{N,\theta}(\bw^n)} \otimes   {D_{N,\theta}(\bw^n)}},   \nabla \mathbb{A}_{\theta}^{}D_{N,\theta}^{}(\bw^n)\right)&=\left(  D_{N,\theta}(\bw^n) \cdot \nabla {D_{N,\theta}(\bw^n)},  D_{N,\theta}^{}(\bw^n)\right) = 0,
\end{array}
\end{equation}
\begin{equation}
 \label{divergencfreebar2}
 \begin{array}{lll}
\left({{D_{N,\theta}(\bw^n)} \otimes  \vec{B}^n} ,   \nabla \vec{B}^n \right)&=\left(  {{D_{N,\theta}(\bw^n)} \cdot \nabla \vec{B}^n},  \vec{B}^n\right) = 0,
\end{array}
\end{equation}

\begin{equation}
 \label{divergencfreebardivergencfreebar}
 \begin{array}{lll}
\left(\overline{  \vec{B}^n \otimes  \vec{B}^n},   \nabla \mathbb{A}_{\theta}^{}D_{N,\theta}^{}(\bw^n)\right)&=\left(  \vec{B}^n \otimes  \vec{B}^n ,  \nabla D_{N,\theta}^{}(\bw^n)\right) = - \left( \vec{B}^n \otimes  {{D_{N,\theta}(\bw^n)} } ,   \nabla \vec{B}^n \right) ,
\end{array}
\end{equation}

 leads to the a priori estimate
\begin{equation}
\begin{array}{lllll}
 \label{apriori1}
\displaystyle \frac{1}{2} \left(\| \mathbb{A}_{\theta}^{\frac{1}{2}}D_{N,\theta}^{\frac{1}{2}}({\bw}^n) \|_{2}^2  + \| {\vec{B}}^n \|_{2}^2 \right)
+ \displaystyle  \int_{0}^{t}  \left(  \nu \| \mathbb{A}_{\theta}^{\frac{1}{2}}D_{N,\theta}^{\frac{1}{2}} {\bw}^n \|_{1,2}^2 +  \mu \|{\vec{B}}^n \|_{1,2}^2\right)  \ ds\\
 \quad  \quad \quad = 
 \displaystyle \frac{1}{2}\left( \| \mathbb{A}_{\theta}^{\frac{1}{2}}D_{N,\theta}^{\frac{1}{2}}\overline{\bv}^{n}_{0} \|_{2}^2   + \| {\vec{B}}^n_0 \|_{2}^2 \right).
\end{array}
\end{equation}
 Using inequality (\ref{lemme13})   we conclude from  (\ref{apriori1}) that 
 \begin{equation}
 \label{iciapriori12}
 \begin{split}
\sup_{t \in [0,T^n[} \left( \| \mathbb{A}_{\theta}^{\frac{1}{2}}D_{N,\theta}^{\frac{1}{2}}({\bw}^n) \|_{2}^2 + \| {\vec{B}}^n \|_{2}^2  \right)+ 2 \int_{0}^{t} \left( \nu  \| \mathbb{A}_{\theta}^{\frac{1}{2}}D_{N,\theta}^{\frac{1}{2}}({\bw}^n) \|_{1,2}^2 + \mu \|{\vec{B}}^n \|_{1,2}^2  \right)\ ds  \le \| {\bv}^{n}_{0} \|_{2}^2 + \| {\vec{B}}^n_0 \|_{2}^2 
\end{split}
\end{equation}
that immediately implies that the existence time is independent of $n$ and it is possible to take $T=T^n$.\\ 
We deduce from (\ref{iciapriori12}) and (\ref{lemme13d}) that 
% \begin{equation}
% \label{vbar1}
%  D_{N,\theta}({\bw}^n) \in L^{\infty}(0,T ; {L}^{2}_{\diver}) \cap L^{2}(0,T ; W^{1,2}(\tore)^3),
%  \end{equation}
%  and 
 \begin{equation}
\label{vbar1sansdn}
\begin{split}
 {\bw}^n \in L^{\infty}(0,T ; \vec{H}^{\theta}_{\sigma }) \cap L^{2}(0,T ; \vec{H}^{1+\theta}_{\sigma }),\\
  \vec{B}^n \in L^{\infty}(0,T ; \vec{L}^{2}_{\sigma }) \cap L^{2}(0,T ; \vec{H}^{1}_{\sigma }).
 \end{split}
 \end{equation}
 Thus it follow from (\ref{lemme13a}) that 
 \begin{equation}
\label{icivbar1}
 D_{N,\theta}({\bw}^n) \in L^{\infty}(0,T ; \vec{H}^{\theta}_{\sigma }) \cap L^{2}(0,T ; \vec{H}^{1+\theta}_{\sigma }),
 \end{equation}

  From  ({\ref{icivbar1}}), (\ref{vbar1sansdn}) and by using H\"{o}lder inequality combined with Sobolev injection we get 
  \begin{equation}
\label{icivtilde1rev}
\begin{split}
{D_{N,\theta}({\bw}^n) \otimes   D_{N,\theta}({\bw}^n)} \in L^{2}(0,T ; H^{2\theta-\frac{1}{2}}( \tore )^{3 \times 3}),\\
{ \vec{B}^n \otimes   \vec{B}^n} \in L^{2}(0,T ; H^{-\frac{1}{2}}( \tore)^{3 \times 3}),\\
{D_{N,\theta}({\bw}^n) \otimes   \vec{B}^n  } \in L^{2}(0,T ; L^{{2}}( \tore )^{3 \times 3}), \textrm{ for any }\theta \ge \frac{1}{2},\\
 { \vec{B}^n  \otimes  D_{N,\theta}({\bw}^n)   } \in L^{2}(0,T ; L^{{2}}( \tore )^{3 \times 3}) \textrm{ for any }\theta \ge \frac{1}{2}.
\end{split} 
 \end{equation}
 From (\ref{icivtilde1rev}) and (\ref{lemme2.2})  it follows   that  
 \begin{equation}
\label{vtilde1}
\begin{split}
\overline{D_{N,\theta}({\bw}^n) \otimes  D_{N,\theta}({\bw}^n)} \in L^{2}(0,T ; H^{4\theta-\frac{1}{2}}( \tore)^{3 \times 3}),\\
\overline{\vec{B}^n \otimes   \vec{B}^n} \in L^{2}(0,T ; H^{2\theta - \frac{1}{2}}( \tore)^{3 \times 3}).
\end{split}  
 \end{equation}
Consequently from the elliptic theory  (\ref{pressuregalerkine}) implies that 
 \begin{equation}
\label{icivbarvbarpressure}
\int_{0}^{T}\|q^n\|_{2\theta-\frac{1}{2},2}^2 dt < K. 
\end{equation}

From  (\ref{weak1galerkine}), (\ref{vbar1sansdn}), (\ref{vtilde1}) and (\ref{icivbarvbarpressure}) we  obtain that 
 \begin{equation}
\label{icivtemps}
\int_{0}^{T}  \|\partial_t\bw^n\|_{2\theta-\frac{3}{2},2}^2  dt < K. 
\end{equation}
Finally, from  (\ref{weakb1galerkine}), (\ref{vbar1sansdn}) and  (\ref{icivtilde1rev}) we also obtain that 
 \begin{equation}
\label{icivbtemps}
\int_{0}^{T}  \|\partial_t\vec{B}^n\|_{-{1},2}^2  dt < K. 
\end{equation}

\textbf{Step 3} (Limit $n \rightarrow \infty$) It follows from the estimates (\ref{vbar1sansdn})-(\ref{icivbtemps}) and the Aubin-Lions compactness lemma
(see \cite{sim87} for example) that there are a  not relabeled  subsequence of $(\bw^n, \vec{B}^n, q^n)$  and a couple $(\bw, \vec{B}, q)$ such that
\begin{align}
\bw^n &\rightharpoonup^* \bw &&\textrm{weakly$^*$ in } L^{\infty}
(0,T;\vec{H}^{\theta}_{\sigma }), \label{c122}\\
\vec{B}^n &\rightharpoonup^* \vec{B} &&\textrm{weakly$^*$ in } L^{\infty}
(0,T;\vec{L}^{{2}}_{\sigma }), \label{bc122}\\
D_{N,\theta}(\bw^n) &\rightharpoonup^* D_{N,\theta}(\bw) &&\textrm{weakly$^*$ in } L^{\infty}
(0,T;\vec{H}^{\theta}_{\sigma }), \label{DNc122}\\
\bw^n &\rightharpoonup \bw &&\textrm{weakly in }
L^2(0,T;\vec{H}^{1+\theta}_{\sigma }), \label{nc22}\\
\vec{B}^n &\rightharpoonup \vec{B} &&\textrm{weakly in }
L^2(0,T;\vec{H}^{{1}}_{\sigma }), \label{bnc22}\\
D_{N,\theta}(\bw^n) &\rightharpoonup D_{N,\theta}(\bw) &&\textrm{weakly in }
L^2(0,T;\vec{H}^{1+\theta}_{\sigma }), \label{DNnc22}\\
\partial_t\bw^n&\rightharpoonup \partial_t\bw &&\textrm{weakly in } L^{2}
(0,T;\vec{H}^{2\theta-\frac{3}{2}}),
\label{nc322}\\
\partial_t\vec{B}^n&\rightharpoonup \partial_t\vec{B} &&\textrm{weakly in } L^{2}
(0,T;\vec{H}^{-{1}}),
\label{bnc322}\\
q^n&\rightharpoonup q &&\textrm{weakly in } L^{2}(0,T;H^{2\theta - \frac{1}{2}}
(\tore)), \label{c32}\\
%\bw^n &\rightharpoonup \bw &&\textrm{weakly in } L^{\frac{8}{3}}
%(0,T;L^{\frac{8}{3}}(\partial \Omega)^3). \label{c5.12}\\
\bw^n &\rightarrow \bw &&\textrm{strongly in  }
L^2(0,T;\vec{H}^{\theta}_{\sigma }),
\label{c83icil2}\\
\vec{B}^n &\rightarrow \vec{B} &&\textrm{strongly in  }
L^2(0,T;\vec{L}^{{2}}_{\sigma }),
\label{bc83icil2}\\
D_{N,\theta}(\bw^n) &\rightarrow D_{N,\theta}(\bw) &&\textrm{strongly in  }
L^2(0,T;\vec{H}^{\theta}_{\sigma }),
\label{DNc83icil2}
%\overline{\bw}^n &\rightarrow \overline{\bw} &&\textrm{strongly in  }
%L^q(0,T;L^q(\tore)^3) \textrm{  for all } q< 5 .\label{c82pppp}
\end{align}
 From (\ref{DNnc22}) and (\ref{DNc83icil2})  it follows   that  
% \begin{align}
%     \overline{D_{N,\theta}({\bw}^n) \times \nabla \times D_{N,\theta}({\bw}^n)} &\rightharpoonup \overline{D_{N,\theta}\bw \times \nabla \times D_{N,\theta} {\bw}} &&\textrm{weakly in }
% L^2(0,T;W^{-\frac{5}{6},2}(\tore)^{3\times 3}), \label{c22prime}\\     
% \overline{D_{N,\theta}({\bw}^n) \times \nabla \times D_{N,\theta} {\bw}^n} &\rightarrow \overline{D_{N,\theta} {\bw} \times \nabla \times D_{N,\theta}\bw} &&\textrm{strongly in  }
% L^1(0,T;L^{1}(\tore)^{3\times 3}),\label{c82int}
%  \end{align}

 \begin{align}
 \overline{D_{N,\theta}({\bw}^n) \otimes  D_{N,\theta}({\bw}^n}) &\rightarrow \overline{D_{N,\theta}({\bw}) \otimes  D_{N,\theta}(\bw)} &&\textrm{strongly in  } L^1(0,T;L^{1}(\tore)^{3 \times 3}),\label{c82int}
 \end{align}

  From (\ref{bnc22}) and (\ref{DNc83icil2})  it follows   that 
  \begin{align}
 \vec{B}^n \otimes   D_{N,\theta}({\bw}^n) &\rightarrow {{\vec{B}} \otimes \bw} &&\textrm{strongly in  }
L^1(0,T;L^{1}(\tore)^{3 \times 3}),\label{bwc82int}\\
 D_{N,\theta}({\bw}^n) \otimes  \vec{B}^n   &\rightarrow \bw \otimes \vec{B}  &&\textrm{strongly in  }
L^1(0,T;L^{1}(\tore)^{3 \times 3}),\label{wbc82int}
 \end{align}

  From (\ref{bnc22}) and (\ref{bc83icil2})  it follows   that 
  \begin{align}
 \overline{\vec{B}^n \otimes   \vec{B}^n} &\rightarrow \overline{\vec{B} \otimes \vec{B}} &&\textrm{strongly in  }
L^1(0,T;L^{1}(\tore)^{3 \times 3}),\label{bc82int}
 \end{align}

 %for all $ q< 2$ and  all $  r < \frac{1}{6}.$\\
Finally, since the
sequence $ \left\{ \overline{D_{N,\theta}({\bw}^n) \otimes D_{N,\theta}({\bw}^n)} \right\}_{n \in \N}$  is bounded in $L^2(0,T;H^{4\theta-\frac{3}{2}}(\tore)^{3 \times 3})$, it converges weakly, up to a
subsequence, to some $\chi $ in  $ L^2(0,T;H^{4\theta-\frac{3}{2}}(\tore)^{3 \times 3})$. The result above (\ref{c82int}) and uniqueness of the limit,
allows us to claim that $\chi  =\overline{D_{N,\theta}(\bw) \otimes  D_{N,\theta}({\bw})}$. Consequently
 \begin{align}
     \overline{D_{N,\theta}({\bw}^n) \otimes D_{N,\theta}({\bw}^n}) &\rightharpoonup \overline{D_{N,\theta}(\bw) \otimes D_{N,\theta}({\bw})} &&\textrm{weakly in }
 L^2(0,T;H^{4\theta-\frac{3}{2}}(\tore)^{3 \times 3}), \label{c22primenc22primen}\    
  \end{align}
  Similarly, 
   \begin{align}
     {D_{N,\theta}({\bw}^n) \otimes \vec{B}^n} &\rightharpoonup {D_{N,\theta}(\bw) \otimes  \vec{B}} &&\textrm{weakly in }
 L^2(0,T;L^{{2}}(\tore)^{3 \times 3}), \label{cb22primencb22primen}\\
  {\vec{B}^n \otimes D_{N,\theta}({\bw}^n)} &\rightharpoonup { \vec{B} \otimes D_{N,\theta}\bw } &&\textrm{weakly in }
 L^2(0,T;L^{{2}}(\tore)^{3 \times 3}), \label{bc22primenbc22primenbc22primen}\\
 \overline{\vec{B}^n \otimes \vec{B}^n } &\rightharpoonup \overline{\vec{B} \otimes \vec{B}} &&\textrm{weakly in }
 L^2(0,T;H^{2\theta-\frac{1}{2}}(\tore)^{3 \times 3}), \label{bc22primenbc22primenbc22primenbc22primen}
  \end{align}

The above established convergences are clearly sufficient for taking the limit in (\ref{weak1galerkine})-(\ref{weakb1galerkine}) and for concluding that   $ \bw$ and $\vec{B}$  satisfy (\ref{weak1})-(\ref{weakb1}). 
Moreover, 
from (\ref{nc22}) and (\ref{nc322})
 one  can deduce by a classical argument (see in \cite{A01})   that 
 \begin{equation}
 \bw \in  \mathcal{C}(0,T;\vec{H}^{\theta}_{\sigma }).
\end{equation}
Similarly, we deduce from (\ref{bnc22}) and (\ref{bnc322}) that 
\begin{equation}
 \vec{B} \in  \mathcal{C}(0,T;\vec{L}^{{2}}_{\sigma }).
\end{equation}
Furthermore, from  the  strong continuity of $\bw$  and $ \vec{B}$ with respect to the time with value in $\vec{H}^{\theta}$ and $\vec{L}^{{2}}_{\sigma }$ respectively,  we deduce   that $\bw(0)=\bw_0$ and $\vec{B}(0)=\vec{B}_{0}$.\\
Let us mention also that  for $\theta \ge \frac{1}{2}$, $D_{N,\theta}({\bw})+ \alpha^{2\theta}(-\Delta)^{\theta}D_{N,\theta}({\bw}) \in L^2(0,T;\vec{H}^{1-\theta}_{\sigma })\hookrightarrow L^2(0,T;\vec{H}^{\frac{3}{2}-2\theta}_{\sigma })$  and $\vec{B} \in L^2(0,T;\vec{H}^{{1}}_{\sigma })$ hence  $\mathbb{A}_{\theta}^{}D_{N,\theta}^{}{\bw} $  is a possible  test function in the weak formulation (\ref{weak1}) and $\vec{B}$ is a possible  test function in the weak formulation (\ref{weakb1}). Thus $ \mathbb{A}_{\theta}^{\frac{1}{2}}D_{N,\theta}^{\frac{1}{2}}({\bw}) $ and $\vec{B}$ verify for all $t \in [0,T]$  the following equality   
\begin{equation}
\begin{array}{lllll}
 \label{apriorileary1}
\displaystyle \frac{1}{2} \left(\| \mathbb{A}_{\theta}^{\frac{1}{2}}D_{N,\theta}^{\frac{1}{2}}({\bw}) \|_{2}^2  + \| {\vec{B}} \|_{2}^2 \right)
+ \displaystyle  \int_{0}^{t}  \left(  \nu \| \mathbb{A}_{\theta}^{\frac{1}{2}}D_{N,\theta}^{\frac{1}{2}}({\bw}) \|_{1,2}^2 +  \mu \|{\vec{B}} \|_{1,2}^2\right)  \ ds\\
 \quad  \quad \quad = 
 \displaystyle \frac{1}{2}\left( \| \mathbb{A}_{\theta}^{\frac{1}{2}}D_{N,\theta}^{\frac{1}{2}}(\overline{\bv}^{}_{0}) \|_{2}^2   + \| {\vec{B}}_0 \|_{2}^2 \right).
\end{array}
\end{equation}

 \textbf{Step 4} (Uniqueness)
%Since the pressure part of the solution is uniquely determined by the velocity part it remain to show the uniqueness to the velocity. 
Next, we will show the continuous dependence of the  solutions on the initial data and in particular the uniqueness.\\
For simplicity we restrict ourselves  to the  critical case $\theta = \frac{1}{2}$, 
 and  we drop some indices of $\theta$. Thus  we will write   ``$D_{N}$'' instead of ``$ D_{N,\theta}$'' and ``$\mathbb{A}_{}$'' instead of ``$ \mathbb{A}_{\theta}$''  expecting that no confusion will occur.\\
Let $({\bw_1,\vec{B}_1},q_1)$ and $({\bw_2,\vec{B}_2},q_2)$   be any two solutions of (\ref{dalpha ns}) on the interval $[0,T]$, with initial values $(\bw_1(0),\vec{B}_1(0))$ and $(\bw_2(0),\vec{B}_2(0))$. Let us denote by  $\delta \vec{w}_{} =\bw_2-\bw_1$, by $\delta \vec{B} = \vec{B}_2 - \vec{B}_1$ and by $\delta q = q_2 - q_1$.

Then one has
\begin{equation}
\label{matin1}
\begin{split}
\partial_{t} \delta\vec{w}- \nu_1
\Delta\delta\vec{w} + \diver(\overline{D_N(\bw_2) \otimes D_N(\bw_2)})- \diver(\overline{D_N(\bw_1) \otimes D_N(\bw_1}))\\
\quad \quad  - \diver(\overline{\vec{B}_2 \otimes \vec{B}_2})+\diver(\overline{\vec{B}_1 \otimes \vec{B}_1})  + \nabla\delta q =0,\\
\partial_{t} \delta\vec{B}  - \nu_2
\Delta\delta\vec{B} + \diver({D_N(\bw_2) \otimes \vec{B}_2})- \diver({D_N(\bw_1) \otimes \vec{B}_1})\\
\quad \quad  - \diver({\vec{B}_2 \otimes D_N(\bw_2)})+\diver({\vec{B}_1 \otimes D_N(\bw_1)})    = 0,
\end{split}
\end{equation}
and $\delta\vec{w} = 0$, $\delta\vec{B} =0$  at initial time. One can take $\alpha^{}(-\Delta)^{\frac{1}{2}}D_N(\delta\vec{w})+ D_N(\delta\vec{w})$
 as test in  the first equation of (\ref{matin1}) and $ \delta\vec{B}$
 as test in the second equations of (\ref{matin1}). Since $
D_N(\vec{w}_1)$ is divergence-free we have 
\begin{equation}
\begin{split}\displaystyle
\int_{0}^{T} \int_{\tore}   D_N(\vec{w}_1) \otimes D_N(\delta \vec{w}) : \nabla D_N(\delta\vec{w})  \hskip 3cm\\
= -\int_{0}^{T} \int_{\tore}   (D_N(\vec{w}_1) \cdot \nabla)D_N(\delta \vec{w}) \cdot D_N(\delta\vec{w}) =0  ,
 \end{split}
\end{equation}
Thus we obtain by using the fact that the averaging operator  commutes with differentiation under periodic boundary conditions 
 \begin{equation}
 \label{soira}
\begin{split}\displaystyle
\int_{0}^{T} \int_{\tore}\left(\diver(\overline{D_N(\bw_2) \otimes D_N(\bw_2)})-\diver(\overline{D_N(\bw_1) \otimes D_N(\bw_1)})\right) \hskip 3cm\\
\qquad  \cdot\left( \alpha^{}(-\Delta)^{\frac{1}{2}}D_N(\delta\vec{w})+ D_N(\delta\vec{w}) \right)\\
= \int_{0}^{T} \int_{\tore}\left(\diver({D_N(\bw_2) \otimes D_N(\bw_2}))-\diver({D_N(\bw_1) \otimes D_N(\bw_1)})\right) \cdot D_N(\delta\vec{w})\\ 
= -\int_{0}^{T} \int_{\tore}  D_N(\delta \vec{w}) \otimes D_N(\vec{w}_2) : \nabla D_N(\delta\vec{w}).
 \end{split}
\end{equation}
Similarly, 
because 
$
D_N(\vec{w}_1)$ is divergence-free we have 
\begin{equation}
\begin{split}\displaystyle
 \int_{0}^{T} \int_{\tore}   D_N(\vec{w}_1) \otimes \delta \vec{B} : \nabla \delta\vec{B}
 = -\int_{0}^{T} \int_{\tore}   (D_N(\vec{w}_1) \cdot \nabla) \delta \vec{B} \cdot \delta\vec{B}  =0  ,
 \end{split}
\end{equation}
and thus we have the following identity 
 \begin{equation}
\begin{split}\displaystyle
\int_{0}^{T} \int_{\tore}\left(\diver({D_N(\bw_2) \otimes \vec{B}_2})-\diver({D_N(\bw_1) \otimes \vec{B}_1})\right)  \cdot \delta\vec{B} \\
= \int_{0}^{T} \int_{\tore}\left(\diver({D_N(\bw_2) \otimes \vec{B}_2})-\diver({D_N(\bw_1) \otimes \vec{B}_1})\right) \cdot\delta\vec{B}\\ 
= -\int_{0}^{T} \int_{\tore} D_N(\delta \vec{w}) \otimes \vec{B}_2 : \nabla \delta\vec{B}.\\
 \end{split}
\end{equation}
Concerning the remaining terms 
we get by integrations by parts and by using  the fact that the averaging operator commutes with differentiation under periodic boundary conditions
\begin{equation}
\begin{split}\displaystyle
\int_{0}^{T} \int_{\tore}\left(-\diver(\overline{\vec{B}_2 \otimes \vec{B}_2})+\diver(\overline{\vec{B}_1 \otimes \vec{B}_1})\right)  \hskip 3cm\\
\qquad  \cdot\left( \alpha^{}(-\Delta)^{\frac{1}{2}}D_N(\delta\vec{w})+ D_N(\delta\vec{w}) \right)  \hskip -1cm\\
= \int_{0}^{T} \int_{\tore}\left(-\diver({\vec{B}_2 \otimes \vec{B}_2})+\diver({\vec{B}_1 \otimes \vec{B}_1})\right) \cdot D_N(\delta\vec{w})\\ 
= \int_{0}^{T} \int_{\tore}   \vec{B}_1 \otimes \delta \vec{B} : \nabla D_N(\delta\vec{w})
 +  \delta \vec{B} \otimes \vec{B}_2 : \nabla D_N(\delta\vec{w}).
 \end{split}
\end{equation}
%= \int_{0}^{T} \int_{\tore}   \vec{W}_1 \otimes \delta \vec{W} : \nabla \delta\vec{w}\\
%\quad + \int_{0}^{T} \int_{\tore}  \delta \vec{W} \otimes \vec{W}_2 : \nabla \delta\vec{w}.\\
and similarly 
\begin{equation}
\label{soirb}
\begin{split}\displaystyle
\int_{0}^{T} \int_{\tore}\left(-\diver({\vec{B}_2 \otimes D_N(\bw_2)})+\diver({\vec{B}_1 \otimes D_N(\bw_1)})\right) \cdot \delta\vec{B} \\
= \int_{0}^{T} \int_{\tore}\left(\diver({\vec{B}_2 \otimes D_N(\bw_2)})-\diver({\vec{B}_1 \otimes D_N(\bw_1)})\right) \cdot\delta\vec{B}\\ 
= \int_{0}^{T} \int_{\tore}  - (\vec{B}_1 \cdot \nabla) D_N(\delta \vec{w}) \cdot  \delta\vec{B} + \delta \vec{B} \otimes D_N(\vec{w}_2) : \nabla \delta\vec{B}.
 \end{split}
\end{equation}
%Therefore, 
%
% \begin{equation}
%\label{soira}\begin{array}{ll}
% \displaystyle {\frac{d }{2dt}} \int_{\tore} | \delta\vec{u}| ^{2} +\displaystyle \nu_1\int_{\tore}| \nabla
%\delta\vec{u} |^{2}  -\displaystyle \int_{\tore} (H_{N} (\B_{1}) \nabla) \delta\B. \delta
%\vec{u}\\
% \hskip 3cm  =\displaystyle -\int_{\tore}(H_{N} (\delta\vec{u}) \nabla) \vec{u}_{2}. \delta
%\vec{u} +\displaystyle \displaystyle\int_{\tore}(H_{N} (\delta\B) \nabla) \B_{2}. \delta
%\vec{u},
%\end{array}
%\end{equation}
%and 
%  \begin{equation}
%\label{soirb}\begin{array}{ll}
% \displaystyle {\frac{d }{2dt}} \int_{\tore}  | \delta\B| ^{2} + \displaystyle\nu_2\int_{\tore} | \nabla
%\delta\B |^{2}   - \displaystyle \int_{\tore} (H_{N} (\B_{1}) \nabla) \delta\vec{u} . \delta
%\B\\ 
% \hskip 3cm  = \displaystyle \int_{\tore}(H_{N} (\delta\B) \nabla) \vec{u}_{2}     . \delta
%\B   \displaystyle   -\displaystyle\int_{\tore}   (H_{N} (\delta\vec{u}) \nabla) \B_{2}  . \delta
%\B.
%\end{array}
%\end{equation}
%
%One has by a integration by parts,
%\BEQ \label{MATIN2}
%- \displaystyle \int_{\tore} (H_{N} (\B_{1}) \nabla) \delta\B. \delta
%\vec{u}     = \displaystyle \int_{\tore} (H_{N} (\B_{1}) \nabla) \delta\vec{u} . \delta
%\B. 
%\EEQ
Therefore by adding (\ref{soira})-(\ref{soirb}) and using  the fact that the averaging operator commutes with differentiation under periodic boundary conditions we obtain
\begin{equation}
 \label{soirab}
\begin{array}{llll}
 \displaystyle \frac{1}{2}\frac{d}{dt} \left(\|\mathbb{A}^{\frac{1}{2}}D_N^{\frac{1}{2}}({\delta \vec{w}})\|_{2}^{2} + \|{\delta \vec{B}}\|_{2}^{2} \right)+\nu \|\nabla \mathbb{A}^{\frac{1}{2}}D_N^{\frac{1}{2}}({\delta \vec{w}})\|_{2}^{2} + \mu \|{\delta \vec{B}}_{}\|_{1,2}^{2} \\
\hskip 1cm=\displaystyle \int_{\tore} D_N( \delta \vec{w}) \otimes D_N(\vec{w}_2) : \nabla D_N(\delta\vec{w})   + \displaystyle \displaystyle \int_{\tore}  D_N(\delta \vec{w}) \otimes \vec{B}_2 : \nabla \delta\vec{B}  \\
\hskip 1.5cm
- \displaystyle \int_{\tore}   \delta \vec{B} \otimes  \vec{B}_2 : \nabla D_N(\delta\vec{w})  - \displaystyle  \int_{\tore}  \delta \vec{B} \otimes  D_N(\vec{w}_2) : \nabla \delta\vec{B}.
 \end{array}
\end{equation}

Next, we estimate the four integrals in the right hand side of (\ref{soirab}). The estimates are obtained by using  H\"{o}lder inequality, Sobolev embedding theorem, Young inequality  and Lemma \ref{lemme1}.

\begin{equation}
 \label{notuniformzero1}
 \begin{split}
   \left| \int_{\tore} D_N( \delta \vec{w}) \otimes D_N(\vec{w}_2) : \nabla D_N(\delta\vec{w}) \right|
  %= \displaystyle \left| \left( D_N(\bw_{2}) - D_N(\bw)_{1} \otimes D_N(\bw)_{1}, \nabla D_N(\delta \vec{w}) \right) \right|\\
 % ( {D_N\bw_{2}\times \nabla \times  {D_N\bw}_{2}}-{D_N\bw_{1}\times \nabla \times {D_N\bw}_{1}}, {D_N\delta \vec{w}} )\\
 \le  \frac{C(N+1)^2}{\alpha\nu} \|{D_N(\delta \vec{w})}\otimes   D_N({\bw}_{2}) \|_{-{\frac{1}{2}},2}^{2} + \frac{\nu \alpha}{4} \|\nabla {\delta \vec{w}} \|_{\frac{1}{2},2}^2 \\
 \hskip 1.5cm \le \displaystyle \frac{C(N+1)^6}{\alpha\nu} \|\delta \vec{w}\|_{\frac{1}{2},2}^{2} \| {\bw}_{2} \|^{2}_{\frac{3}{2},2} + \frac{\nu \alpha}{4} \|\nabla {\delta \vec{w}} \|_{\frac{1}{2},2}^2,
% + \displaystyle \frac{C(\alpha)(N+1)^2}{\nu} \|{ \delta \vec{B}}\otimes   {\vec{B}}_{2} \|_{-{\frac{1}{2}},2}^{2} + \frac{\nu \alpha}{4} \|\nabla {\delta \vec{w}} \|_{\frac{1}{2},2}^2 \\
% + \displaystyle \frac{C}{\mu} \|{D_N \delta \vec{w}}\otimes   {\vec{B}}_{1} \|_{2}^{2} + \frac{\mu}{4} \|\nabla {\delta \vec{B}} \|_{2}^2 \\
% + \displaystyle \frac{C}{\mu} \| \delta {\vec{B}}_{} \otimes {D_N  \vec{w}_{1}}   \|_{2}^{2} + \frac{\mu }{4} \|\nabla {\delta \vec{B}} \|_{2}^2. 
\end{split}
\end{equation}
\begin{equation}
 \label{notuniformzero2}
 \begin{split}
  \displaystyle \left| \int_{\tore}  D_N(\delta \vec{w}) \otimes \vec{B}_2 : \nabla \delta\vec{B}   \right| 
  %= \displaystyle \left| \left( D_N(\bw_{2}) - D_N(\bw)_{1} \otimes D_N(\bw)_{1}, \nabla D_N(\delta \vec{w}) \right) \right|\\
 % ( {D_N\bw_{2}\times \nabla \times  {D_N\bw}_{2}}-{D_N\bw_{1}\times \nabla \times {D_N\bw}_{1}}, {D_N\delta \vec{w}} )\\
% \le \displaystyle \frac{C(\alpha)(N+1)^2}{\nu} \|{D_N \delta \vec{w}}\otimes   {D_N \bw}_{1} \|_{-{\frac{1}{2}},2}^{2} + \frac{\nu \alpha}{4} \|\nabla {\delta \vec{w}} \|_{\frac{1}{2},2}^2 \\
% + \displaystyle \frac{C(\alpha)(N+1)^2}{\nu} \|{ \delta \vec{B}}\otimes   {\vec{B}}_{1} \|_{-{\frac{1}{2}},2}^{2} + \frac{\nu \alpha}{4} \|\nabla {\delta \vec{w}} \|_{\frac{1}{2},2}^2 \\
% \le  \displaystyle \frac{C}{\mu} \| \delta {\vec{B}}_{} \otimes {D_N  \vec{w}_{1}}   \|_{2}^{2} + \frac{\mu }{4} \|\nabla {\delta \vec{B}} \|_{2}^2. 
\le  \displaystyle \frac{C}{\mu} \|{D_N(\delta \vec{w})}\otimes   {\vec{B}}_{2} \|_{2}^{2} + \frac{\mu}{4} \|\nabla {\delta \vec{B}} \|_{2}^2\\
\hskip 1.5cm  \le  \displaystyle \frac{C(N+1)^2}{\mu} \| \delta \vec{w}\|_{\frac{1}{2},2}^2 \|{\vec{B}}_{2} \|_{1,2}^{2} + \frac{\mu}{4} \|\nabla {\delta \vec{B}} \|_{2}^2, 
\end{split}
\end{equation}

\begin{equation}
 \label{notuniformzero3}
 \begin{split}
  \displaystyle \left|  \int_{\tore}   \delta \vec{B} \otimes  \vec{B}_2 : \nabla D_N(\delta\vec{w}) \right|  
  %= \displaystyle \left| \left( D_N(\bw_{2}) - D_N(\bw)_{1} \otimes D_N(\bw)_{1}, \nabla D_N(\delta \vec{w}) \right) \right|\\
 % ( {D_N\bw_{2}\times \nabla \times  {D_N\bw}_{2}}-{D_N\bw_{1}\times \nabla \times {D_N\bw}_{1}}, {D_N\delta \vec{w}} )\\
\le \displaystyle \frac{C(N+1)^2}{\alpha \nu} \|    \delta \vec{B} \otimes  \vec{B}_2 \|_{-{\frac{1}{2}},2}^{2} + \frac{\nu \alpha}{4} \|\nabla {\delta \vec{w}} \|_{\frac{1}{2},2}^2 \\
 \hskip 1.5cm \le \displaystyle \frac{C(N+1)^2}{\alpha \nu} \|\delta \vec{B}\|_{2}^{2} \| \vec{B}_{2} \|^{2}_{1,2} + \frac{\nu \alpha}{4} \|\nabla {\delta \vec{w}} \|_{\frac{1}{2},2}^2,
% + \displaystyle \frac{C(\alpha)(N+1)^2}{\nu} \|{ \delta \vec{B}}\otimes   {\vec{B}}_{2} \|_{-{\frac{1}{2}},2}^{2} + \frac{\nu \alpha}{4} \|\nabla {\delta \vec{w}} \|_{\frac{1}{2},2}^2 \\
% + \displaystyle \frac{C}{\mu} \|{D_N \delta \vec{w}}\otimes   {\vec{B}}_{1} \|_{2}^{2} + \frac{\mu}{4} \|\nabla {\delta \vec{B}} \|_{2}^2 \\
% + \displaystyle \frac{C}{\mu} \| \delta {\vec{B}}_{} \otimes {D_N  \vec{w}_{1}}   \|_{2}^{2} + \frac{\mu }{4} \|\nabla {\delta \vec{B}} \|_{2}^2. 
\end{split}
\end{equation}

\begin{equation}
 \label{notuniformzero4}
 \begin{split}
  \displaystyle \left|    \int_{\tore}  \delta \vec{B} \otimes  D_N(\vec{w}_2) : \nabla \delta\vec{B}  \right|=  \displaystyle \left| \int_{\tore}    (\delta \vec{B} \cdot \nabla ) D_N(\vec{w}_2) \cdot  \delta\vec{B}  \right| \\
  %= \displaystyle \left| \left( D_N(\bw_{2}) - D_N(\bw)_{1} \otimes D_N(\bw)_{1}, \nabla D_N(\delta \vec{w}) \right) \right|\\
 % ( {D_N\bw_{2}\times \nabla \times  {D_N\bw}_{2}}-{D_N\bw_{1}\times \nabla \times {D_N\bw}_{1}}, {D_N\delta \vec{w}} )\\
% \le \displaystyle \frac{C(\alpha)(N+1)^2}{\nu} \|{D_N \delta \vec{w}}\otimes   {D_N \bw}_{1} \|_{-{\frac{1}{2}},2}^{2} + \frac{\nu \alpha}{4} \|\nabla {\delta \vec{w}} \|_{\frac{1}{2},2}^2 \\
% + \displaystyle \frac{C(\alpha)(N+1)^2}{\nu} \|{ \delta \vec{B}}\otimes   {\vec{B}}_{1} \|_{-{\frac{1}{2}},2}^{2} + \frac{\nu \alpha}{4} \|\nabla {\delta \vec{w}} \|_{\frac{1}{2},2}^2 \\
% \le  \displaystyle \frac{C}{\mu} \| \delta {\vec{B}}_{} \otimes {D_N  \vec{w}_{1}}   \|_{2}^{2} + \frac{\mu }{4} \|\nabla {\delta \vec{B}} \|_{2}^2. 
\hskip 1.5cm \le \displaystyle \frac{C}{\mu} \|{\delta \vec{B}}\|_{L^2}^2 \| \nabla D_N(\vec{w}_2) \|_{L^3}^{2} + \frac{\mu}{4} \| {\delta \vec{B}} \|_{L^6}^2\\
\hskip 1.5cm  \le  \displaystyle \frac{C(N+1)^2}{\mu} \|{\delta \vec{B}}\|_{2}^2 \|  \vec{w}_2 \|_{\frac{3}{2},2}^{2} + \frac{\mu}{4} \| {\delta \vec{B}} \|_{1,2}^2.
\end{split}
\end{equation}

By using  (\ref{lemme13d}) we have 
\begin{equation}
\label{notuniform}
 \begin{split}
  \displaystyle 
   \frac{1}{2}\frac{d}{dt}\left( \alpha^{} \|{\delta \vec{w}}\|_{\frac{1}{2},2}^{2} + \|{\delta \vec{B}}\|_{2}^{2} \right)+ \nu \alpha^{} \|\nabla {\delta \vec{w}}_{}\|_{\frac{1}{2}}^{2}  + \mu \|  {\delta \vec{B}}_{}\|_{1,2}^{2}  \hskip 3cm\\
\le 
  \frac{1}{2}\frac{d}{dt} \left( \|\mathbb{A}^{\frac{1}{2}}D_N^{\frac{1}{2}}({\delta \vec{w}})\|_{2}^{2}+ \|{\delta \vec{B}}\|_{2}^{2}\right) +\nu \|\nabla \mathbb{A}^{\frac{1}{2}}D_N^{\frac{1}{2}}({\delta \vec{w}})\|_{2}^{2}  + \mu \|  {\delta \vec{B}}_{}\|_{1,2}^{2} \end{split}
\end{equation}
From (\ref{soirab})-(\ref{notuniform})  we get 
\begin{equation}
 \begin{split}
  \displaystyle 
  \frac{d}{dt}\left(\alpha^{}  \|{\delta \vec{w}}\|_{\frac{1}{2},2}^{2} + \|{\delta \vec{B}}\|_{2}^{2} \right) + \nu \alpha^{} \|\nabla {\delta \vec{w}}_{}\|_{\frac{1}{2}}^{2}   + \mu \|{\delta \vec{B}}_{}\|_{1	,2}^{2}  \hskip 3cm \\
 \hskip 1cm \le \displaystyle \frac{C(N+1)^6}{\alpha \nu} \|\delta \vec{w}\|_{\frac{1}{2},2}^{2} \| {\bw}_{2} \|^{2}_{\frac{3}{2},2} +
  \displaystyle \frac{C(N+1)^2}{\mu} \| \delta \vec{w}\|_{\frac{1}{2},2}^2 \|{\vec{B}}_{2} \|_{1,2}^{2} \\
 \hskip 1cm +\frac{C(N+1)^2}{\alpha \nu} \|\delta \vec{B}\|_{2}^{2} \| \vec{B}_{2} \|^{2}_{1,2} 
+  \displaystyle \frac{C(N+1)^2}{\mu} \|{\delta \vec{B}}\|_{2}^2 \|  \vec{w}_2 \|_{\frac{3}{2},2}^{2}.
\end{split}
\end{equation}

Hence, 

\begin{equation}
 \begin{split}
  \displaystyle 
  \frac{d}{dt}\left(\alpha^{}  \|{\delta \vec{w}}\|_{\frac{1}{2},2}^{2} + \|{\delta \vec{B}}\|_{2}^{2} \right) + \nu \alpha^{} \|\nabla {\delta \vec{w}}_{}\|_{\frac{1}{2}}^{2}   + \mu \|{\delta \vec{B}}_{}\|_{1	,2}^{2}  \hskip 3cm \\
 \hskip 1cm \le \displaystyle \frac{C(N+1)^6}{\min{(\alpha \nu, \mu)}} \left( \alpha \|\delta \vec{w}\|_{\frac{1}{2},2}^{2} + \|\delta \vec{B}\|_{2}^{2} \right) \left( \| {\bw}_{2} \|^{2}_{\frac{3}{2},2} +
  \displaystyle   \|{\vec{B}}_{2} \|_{1,2}^{2} \right).
\end{split}
\end{equation}

 % \begin{equation}
% \begin{split}
%  \displaystyle 
%\|{\delta \vec{w}}_{,t}\|_{2}^{2} +\nu \|\nabla{\delta \vec{w}}_{}\|_{2}^{2} & \le \displaystyle
%\frac{4}{\nu} \|\widetilde{\delta \vec{w}}
%\overline{\bw}_{1} \|_{2}^{2} +  \frac{4}{\nu} \|\widetilde{\vec{v}_2}
%\overline{\delta \vec{w}} \|_{2}^{2}    \\
%&\le \displaystyle \frac{4}{\nu} \|\widetilde{\delta \vec{w}}\|_{\frac{1}{2}-2\theta_2,2}^{2}
%\|\overline{\bw}_{1} \|^{2}_{1+2\theta_2}+  \displaystyle\frac{4}{\nu} \|\overline{\delta \vec{w}}\|_{\frac{1}{2}-2\theta_1,2}^{2}
%\|\widetilde{\bw}_{2} \|^{2}_{1+2\theta_1}  \\ &\le \displaystyle \frac{1}{\alpha^{{2\theta_1}+{2\theta_2}}}\frac{4}{\nu} \|{\delta \vec{w}}\|_{2}^{2}
%\left(\|{\bw}_{1} \|^{2}_{1,2} + 
%\|{\bw}_{2} \|^{2}_{1,2}\right)
%\end{split}
%\end{equation}
%\begin{equation}
% \begin{split}
%  \displaystyle 
%\|{\delta \vec{w}}_{,t}\|_{2}^{2} +\nu \|\nabla{\delta \vec{w}}_{}\|_{2}^{2} & \le \displaystyle
%\frac{4}{\nu} \|\widetilde{\delta \vec{w}}
%\overline{\bw}_{1} \|_{2}^{2} +  \frac{4}{\nu} \|\widetilde{\vec{v}_2}
%\overline{\delta \vec{w}} \|_{2}^{2}    \\
%&\le \displaystyle \frac{4}{\nu} \|\widetilde{\delta \vec{w}}\|_{\frac{3}{4}+2\theta_1-\theta_2,2}^{2}
%\|\overline{\bw}_{1} \|^{2}_{\frac{3}{4}-2\theta_1+\theta_2}+  \displaystyle\frac{4}{\nu} \|\overline{\delta \vec{w}}\|_{\frac{1}{2}-2\theta_1+\theta_2,2}^{2}
%\|\widetilde{\bw}_{2} \|^{2}_{1+2\theta_1-\theta_2}  \\ &\le \displaystyle \frac{1}{\alpha^{{2\theta_1}+{2\theta_2}}}\frac{4}{\nu} \|{\delta \vec{w}}\|_{2}^{2}
%\left(\|{\bw}_{1} \|^{2}_{1,2} + 
%\|{\bw}_{2} \|^{2}_{1,2}\right)
%\end{split}
%\end{equation}
 Since $\| {\bw}_{2} \|^{2}_{\frac{3}{2},2} +\|{\vec{B}}_{2} \|_{1,2}^{2} \in L^1([0,T])$, we get by using Gronwall's inequality  
the continuous dependence of the solutions on the initial data. In particular, if ${\delta \vec{w}}^{}_{0}= {\delta \vec{B}}^{}_{0}=0$ then ${\delta \vec{w}}= {\delta \vec{B}}=0$ and the solutions are unique for all $t \in [0,T] .$  Since $T>0$ is arbitrary this solution may be uniquely extended for all time.\\ 
This finishes the proof of Theorem \ref{TH1}.\\

\section{Regular weak solution in the  inviscid case ($\nu=0$, $\mu >0$)}
\label{sec4}

In this section we prove  Theorem \ref{Theorem2} by using the  Galerkin method. 

\textbf{Step 1:} (Existence of weak solutions)

\textbf{Step:1-i}(Galerkin approximation).
Consider the  sequence $\left\{ \bfi^{r} \right\}_{r=1}^{\infty}$ defined in the proof of Theorem \ref{TH1}.  
% We set 
% \begin{equation}
% \begin{split}
% \bw^n(t,\bx)=\sum_{{r=1}}^{n}\bc_{r}^n (t) \bfi^{r}(\bx),\\
% \vec{B}^n(t,\bx)=\sum_{{r=1}}^{n}\vec{d}_{r}^n (t) \bfi^{r}(\bx),
% \\
% \quad  \hbox{ and } q^n(t,\bx)=\sum_{{|\vec{k}|=1}}^{n}q_{\vec{k}}^n (t) e^{i \vec{k} \cdot \bx}.
% \end{split}
% \end{equation}
% such that $\bk \cdot  \bc_{r}^n = 0  $ for all $\bk \in \Z^3\setminus\{0\}$ and $\left(\bc_{r}^n \right)^{*}=\bc_{-r}^n  $, where  $ \left(\bc_{r}^n \right)^{*}$ denote the complex conjugate $\bc_{r}^n.  $
% Thus due of (\ref{TKEbis}) and (\ref{TKE}) we have 
%\begin{equation}
% \widetilde{\bw}^n(t,\bx)=\sum_{{|r|=1}}^{N}\widetilde{\bc}_{r}^n (t) \bfi^{r}(\bx) \textrm{ and }  \overline{\bw}^n(t,\bx)=\sum_{{|r|=1}}^{N}\overline{\bc}_{r}^n (t) \bfi^{r}(\bx), 
%\end{equation}
% where 
%\begin{equation}
% \widetilde{\bc}_{r}^n = \frac{{\bc}_{r}^n }{1+\alpha^{2\theta_1}|r|^{2\theta_1}} \textrm{ and }  \overline{\bc}_{r}^n = \frac{{\bc}_{r}^n }{1+\alpha^{2\theta_2}|r|^{2\theta_2}}, 
%\end{equation}
%for all $\bk \in \Z^3\setminus\{0\}$.\\
We look for $(\bw^n(t,\bx),\vec{B}^n(t,\bx), q^n(t,\bx)) $, where  \begin{equation}
\begin{split}
\bw^n(t,\bx)=\sum_{{r=1}}^{n}\bc_{r}^n (t) \bfi^{r}(\bx),\\
\vec{B}^n(t,\bx)=\sum_{{r=1}}^{n}\vec{d}_{r}^n (t) \bfi^{r}(\bx),
\\
\quad  \hbox{ and } q^n(t,\bx)=\sum_{{|\vec{k}|=1}}^{n}q_{\vec{k}}^n (t) e^{i \vec{k} \cdot \bx},
\end{split}
\end{equation} 
 that are determined through the system of equations 
\begin{equation}
\begin{split}
 \left( \partial_t\bw^n, \bfi^{r} \right)  -  \left(\overline{D_{N,\theta}(\bw^n) \otimes  D_{N,\theta}(\bw^n)} , 
\nabla\bfi^{r}\right) \; \\
+ \left( \overline{\vec{B}^n \otimes  \vec{B}^n} , \nabla \bfi^{r} \right) = 0 \; ,
\quad {r=1},2,...,n,
\end{split}\label{nuw1g}
\end{equation}
\begin{equation}
\begin{split}
 \left( \partial_t\vec{B}^n, \bfi^{r} \right)  +  \left({D_{N,\theta}(\bw^n) \otimes  \vec{B}^n} , 
\nabla\bfi^{r}\right) +  \mu \left(
\nabla \vec{B}^n, \nabla \bfi^{r} \right)\; \\
- \left( {\vec{B}^n \otimes  D_{N,\theta}(\bw^n)} , \nabla \bfi^{r} \right) = 0 \; ,
\quad {r=1},2,...,n,
\end{split}\label{nuweakb1galerkine}
\end{equation}

and 
%\begin{equation}
%\begin{split}
% q^n = - \sum_{i,j}\partial_i\partial_j \Delta^{-1}(\Pi^n(\overline{{v}_i^n {v}_j^n)})= - \sum_{i,j} R_{ij}(\Pi^n(\widetilde{v}_i^n\overline{v}_j^n)),
%\end{split}\label{pressuregalerkine}
%\end{equation}
\begin{equation}
\begin{split}
 %\displaystyle \Delta^{} q^n=  \nabla \cdot   \Pi^n \left(  \overline{D_{N} {\bw}^n  \times \nabla \times D_{N} {\bw}^n }  \right).
\displaystyle \Delta^{} q^n=  -\diver \diver  \left(\Pi^n(\overline{D_{N,\theta}({\bw}^n) \otimes D_{N,\theta}({\bw}^n)} - \overline{\vec{B}^n \otimes  \vec{B}^n})  \right).
\end{split}\label{nupgnu}
\end{equation}
Where the projector $ \displaystyle \Pi^n $ assign to any Fourier series $\displaystyle \sum_{\bk \in \Z^3\setminus\{0\}} \vec{g}_{\bk} e^{i\bk \cdot\bx} $ the following  series  
 $\displaystyle \sum_{\bk \in \Z^3\setminus\{0\}, |\bk| \le n} \vec{g}_{\bk} e^{i\bk \cdot\bx}. $
 %we deduce from the elliptic equation (\ref{pressuregalerkine}) that 
 % after dropping indices of $N$ for shortness, that 
  % \begin{equation}
%\begin{split}
% p^n = - \sum_{i,j}\partial_i\partial_j \Delta^{-1}(\Pi^n(\widetilde{v}_i^n\overline{v}_j^n))= - \sum_{i,j} R_{ij}(\Pi^n(\widetilde{v}_i^n\overline{v}_j^n)),
%\end{split}\label{pressurepseudo}
%\end{equation}
% and $R_{ij}$ is the Riez operator defined through the Fourier transform by
% \begin{equation}
%\begin{split}
%  \widehat{R_{ij}(u)}=\frac{k_i k_j}{|\bk|^2}\widehat{u({\bk})}, \quad \hbox{ for all } \bk \in \Z^3\setminus\{0\}.
%  \end{split}\label{pressurepseudo2}
%\end{equation}

% The Laplace operator $\displaystyle \Delta^{-1}$ can be  defined through the Fourier transform hence $\displaystyle (\Delta^{-1}\diver \diver)$ may be viewed as a pseudofifferentiel operator 
%  \begin{equation}
%\begin{split}
%   ((\Delta^{-1}\diver \diver){\bw})_{\bk}=\frac{\bk^2}{|\bk|^2}\bw_{\bk}.
%\end{split}\label{pressurepseudo}
%\end{equation}
Moreover we require that $\bw^n $ and $\vec{B}^n $  satisfy the following initial conditions
\begin{equation}
\label{nuinitial Galerkine}
\bw^n(0,.)= \bw^n_0= \sum_{r=1}^{n}\bc_0^n  \bfi^{r}(\bx), \quad  \vec{B}^n(0,.)= \vec{B}^n_0= \sum_{r=1}^{n}\vec{d}_0^n  \bfi^{r}(\bx)
\end{equation}
and 
\begin{equation}
\begin{split}
 \label{nuinitial2 Galerkine}
\bw^n_0 \rightarrow \bw_0  \quad \textrm{ strongly  in }  \vec{H}^{2\theta}_{\sigma} \quad \textrm{ when } n \rightarrow \infty,\\
\vec{B}^n_0 \rightarrow \vec{B}_0  \quad \textrm{ strongly  in }  \vec{L}^{2}_{\sigma} \quad \textrm{ when } n \rightarrow \infty.
\end{split}
\end{equation}
%Where the initial condition $ \overline{\bw}^n_0$ is deduced from $\bw^n_0$ through the relation (\ref{TKE}).
%\begin{equation}
% \label{initialbar Galerkine}
%  \overline{\bw}^n_0   + \alpha^{2\theta}(-\Delta)^{\theta} \overline{\bw}^n_0=\bw^n_0
%\end{equation}
 The classical Caratheodory theory \cite{Wa70}  then implies the short-time existence of solutions 
to (\ref{nuw1g})-(\ref{nupgnu}).  Next we derive  estimates on $\bc^n$ and $\vec{d}^n$ that are uniform w.r.t. $n$.
These estimates then imply that the  solution of  (\ref{nuw1g})-(\ref{nupgnu}) constructed on a short time interval $[0, T^n[ $ exists for all $t \in [0, T]$.\\

\textbf{Step 1-ii:}(Uniform estimates) We are going to derive some new uniform estimates on $(\vec{w}^{n}, \vec{B}^{n}, q^{n})$. 

Multiplying the $r$th equation in (\ref{nuw1g}) with $\alpha^{2 \theta}|\bk|^{2 \theta }\widehat{D_N}{\bc}^n_{r}(t)+\widehat{D_N}{\bc}^n_{r}(t)$,  and the $r$th equation in (\ref{nuweakb1galerkine}) with ${\vec{d}}^n_{r}(t)$ summing over ${r=1},2,...,n$, integrating over time from $0$ to $t$
  we obtain in a similar way than  in the proof of Theorem \ref{TH1},  
 \begin{equation}
 \label{iciapriori12nu}
 \begin{split}
\sup_{t \in [0,T]} \left( \| \mathbb{A}^{\frac{1}{2}}_{\theta}D_{N,\theta}^{\frac{1}{2}}({\bw}^{n}) \|_{2}^2 + \| {\vec{B}}^{n} \|_{2}^2  \right)+ 2 \mu \int_{0}^{t}   \|{\vec{B}}^{n} \|_{1,2}^2  \ ds  \le \| {\bv}^{}_{0} \|_{2}^2 + \| {\vec{B}}_0 \|_{2}^2. 
\end{split}
\end{equation}

The above inequality  implies that the existence time is independent of $n$ and it is possible to take $T=T^n$.\\ 
As a consequence from (\ref{iciapriori12nu}) and Lemma (\ref{lemme13d}), we derive that there exists $C>0$  such that 
 \begin{align}
 \label{nuuniform}
    \| D_N({\bw}^{n}) \|_{L^{\infty}(0,T;\vec{H}^{{\theta}}_{\sigma})}+ \| {\bw}^{n} \|_{L^{\infty}(0,T;\vec{H}^{{\theta}}_{\sigma})} +   \| {\vec{B}}^{n} \|_{L^{\infty}(0,T;\vec{L}^{2}_{\sigma})}   \le C \\
    \label{nuuniform2}
    \| {\vec{B}}^{n} \|_{L^{2}(0,T;\vec{H}^{1}_{\sigma})} \le C
\end{align}

From  ({\ref{nuuniform}})-({\ref{nuuniform2}}) , and by using H\"{o}lder inequality combined with Sobolev injection we get 
  \begin{equation}
\label{nuicivtilde1rev}
\begin{split}
{D_{N,\theta}({\bw}^{n}) \otimes   D_{N,\theta}({\bw}^{n}}) \in L^{\infty }(0,T ; H^{ 2 \theta - \frac{3}{2}}( \tore )^{3 \times 3}),\\
{ \vec{B}^{n} \otimes   \vec{B}^{n}} \in L^{4}(0,T ; H^{-{1}}( \tore)^{3 \times 3})\cap  L^{\frac{6}{5}}(0,T ; H^{\frac{1}{6}}(\tore)^{3 \times 3}),\\
{D_{N,\theta}({\bw}^{n}) \otimes   \vec{B}^{n}  } \in L^{2}(0,T ; H^{\theta - \frac{1}{2} }( \tore )^{3 \times 3}),\\
 { \vec{B}^{n} \otimes  D_{N,\theta}({\bw}^{n})    } \in L^{2}(0,T ; H^{\theta - \frac{1}{2} }( \tore )^{3 \times 3}).
\end{split}
\end{equation} 
Note that $  L^{2}(0,T ; H^{\theta - \frac{1}{2} }( \tore )^{3 \times 3}) \hookrightarrow L^{2}(0,T ; L^{{2} }( \tore )^{3 \times 3}) $  when $ \theta \ge \frac{1}{2}$ thus 
\begin{equation}
\label{nu2icivtilde1rev}
\begin{split}
{D_{N,\theta}({\bw}^{n}) \otimes   \vec{B}^{n}  } \in L^{2}(0,T ; L^{{2} }( \tore )^{3 \times 3}),\\
 { \vec{B}^{n} \otimes  D_{N,\theta}({\bw}^{n}}) \in L^{2}(0,T ; L^{ {2} }( \tore )^{3 \times 3}).
\end{split}
 \end{equation}
 From (\ref{nuicivtilde1rev}) and (\ref{lemme2.2})  it follows   that, for any   $\theta \ge  \frac{5}{6}$ 
 \begin{equation}
\label{nuvtilde1}
\begin{split}
\overline{D_{N,\theta}({\bw}^{n}) \otimes  D_{N,\theta}({\bw}^{n})} &\in L^{\infty }(0,T ; H^{ 4 \theta - \frac{3}{2}}_{}) \hookrightarrow  L^{4}(0,T ; H^{2 \theta - 1}( \tore)^{3 \times 3})\cap  L^{\frac{6}{5}}(0,T ; H^{2 \theta + \frac{1}{6}}( \tore)^{3 \times 3}),\\
%\overline{\vec{B}^n \otimes   \vec{B}^n} \in L^{2}(0,T ; H^{2 \theta - \frac{1}{2}}( \tore)^{3 \times 3}).
\overline{\vec{B}^{n}\otimes   \vec{B}^{n}} &\in L^{4}(0,T ; H^{2 \theta - 1}( \tore)^{3 \times 3})\cap  L^{\frac{6}{5}}(0,T ; H^{2 \theta + \frac{1}{6}}( \tore)^{3 \times 3}).
\end{split}  
 \end{equation}

As in the proof of Theorem \ref{TH1} and by using (\ref{nuuniform})-(\ref{nuicivbarvbarpressure}),  we can derive the following estimates
 \begin{align}
 \label{nuicivbarvbarpressure}
\int_{0}^{T}\|q^{n}\|_{2 \theta - 1 ,2}^4 dt < C, \\
\label{nuicivtemps}
\int_{0}^{T}  \|\partial_t{\bw^{n}}\|_{2 \theta - 2,2}^4  dt < C \textrm{ and }  \int_{0}^{T}  \|\partial_t{\bw^{n}}\|_{2 \theta - \frac{5}{6}}^{\frac{6}{5}}  dt < C, \\
\label{nuicivbtemps}
\int_{0}^{T}  \|\partial_t{\vec{B}^{n}}\|_{-1,2}^2  dt < C. 
\end{align}

\textbf{Step 1-iii} (Limit $n \rightarrow \infty$ part I) It follows from the estimates (\ref{nuuniform})-(\ref{nuicivbtemps}) and the Aubin-Lions compactness lemma
(see \cite{sim87} for example) that there are a  not relabeled  subsequence of $(\bw^{n}, \vec{B}^{n}, q^{n},)$  and a triplet $(\bw, \vec{B}, q)$ such that
\begin{align}
\bw^{n} &\rightharpoonup^* \bw &&\textrm{weakly$^*$ in } L^{\infty}
(0,T;\vec{H}^{\theta}_{\sigma }), \label{nuc122}\\
\vec{B}^{n} &\rightharpoonup^* \vec{B} &&\textrm{weakly$^*$ in } L^{\infty}
(0,T;\vec{L}^{{2}}_{\sigma }), \label{nubc122}\\
D_{N,\theta} (\bw^{n}) &\rightharpoonup^* D_{N,\theta}(\bw) &&\textrm{weakly$^*$ in } L^{\infty}
(0,T;\vec{H}^{\theta}_{\sigma }), \label{nuDNc122}\\
% \bw^n &\rightharpoonup \bw &&\textrm{weakly in }
% L^2(0,T;\vec{H}^{\frac{3}{2}}_{\sigma }), \label{nc22}\\
\vec{B}^{n} &\rightharpoonup \vec{B} &&\textrm{weakly in }
L^2(0,T;\vec{H}^{{1}}_{\sigma }), \label{nubnc22}\\
% D_N\bw^n &\rightharpoonup D_N\bw &&\textrm{weakly in }
% L^2(0,T;\vec{H}^{\frac{3}{2}}_{\sigma }), \label{DNnc22}\\
\partial_t{\bw^{n}}&\rightharpoonup \partial_t\bw &&\textrm{weakly in } L^{4}
(0,T;\vec{H}^{2\theta -2}) \cap L^{\frac{6}{5}}(0,T ; \vec{H}^{2 \theta - \frac{5}{6}}),
\label{nunc322}\\
\partial_t{\vec{B}^{n}}&\rightharpoonup \partial_t\vec{B} &&\textrm{weakly in } L^{2}
(0,T;\vec{H}^{-{1}}),
\label{nubnc322}\\
q^{n}&\rightharpoonup q &&\textrm{weakly in } L^{4}(0,T;H^{2\theta -1}
(\tore)), \label{nuc32}\\
%\bw^n &\rightharpoonup \bw &&\textrm{weakly in } L^{\frac{8}{3}}
%(0,T;L^{\frac{8}{3}}(\partial \Omega)^3). \label{c5.12}\\
\bw^{n} &\rightarrow \bw &&\textrm{strongly in  }
L^2(0,T;\vec{L}^{{2}}_{\sigma }),
\label{nuc83icil2}\\
\vec{B}^{n} &\rightarrow \vec{B} &&\textrm{strongly in  }
L^2(0,T;\vec{L}^{{2}}_{\sigma }),
\label{nubc83icil2}\\
D_{N,\theta}(\bw^{n}) &\rightarrow D_{N,\theta}(\bw) &&\textrm{strongly in  }
L^2(0,T;\vec{L}^{2}_{\sigma }),
\label{nuDNc83icil2}
%\overline{\bw}^n &\rightarrow \overline{\bw} &&\textrm{strongly in  }
%L^q(0,T;L^q(\tore)^3) \textrm{  for all } q< 5 .\label{c82pppp}
\end{align}

 From  (\ref{nuDNc83icil2})  it follows   that

 \begin{align}
 \overline{D_{N,\theta}({\bw}^{n}) \otimes  D_{N,\theta}({\bw}^{n})} &\rightarrow \overline{D_{N,\theta}({\bw}) \otimes  D_{N,\theta}(\bw)} &&\textrm{strongly in  } L^1(0,T;L^{1}(\tore)^{3 \times 3}),\label{nuc82int}
 \end{align}

  From (\ref{nubnc22}) and (\ref{nuDNc83icil2})  it follows   that 
  \begin{align}
 \vec{B}^{n} \otimes   D_{N,\theta}({\bw}^{n}) &\rightarrow {{\vec{B}} \otimes D_{N,\theta}(\bw)} &&\textrm{strongly in  }
L^1(0,T;L^{1}(\tore)^{3 \times 3}),\label{nubwc82int}\\
 D_{N,\theta}({\bw}^{n}) \otimes  \vec{B}^{n}   &\rightarrow D_{N,\theta}(\bw) \otimes \vec{B}  &&\textrm{strongly in  }
L^1(0,T;L^{1}(\tore)^{3 \times 3}),\label{nuwbc82int}
 \end{align}

  From (\ref{nubnc22}) and (\ref{nubc83icil2})  it follows   that 
  \begin{align}
 \overline{\vec{B}^{n} \otimes   \vec{B}^{n}} &\rightarrow \overline{\vec{B} \otimes \vec{B}} &&\textrm{strongly in  }
L^1(0,T;L^{1}(\tore)^{3 \times 3}),\label{bc82intbc82int}
 \end{align}

 %for all $ q< 2$ and  all $  r < \frac{1}{6}.$\\
 Since the
sequence $  \left\{\overline{D_{N,\theta}({\bw}^{n})\otimes D_{N,\theta}({\bw}^{n})} \right\}_{n \in \N}$  is bounded in $L^4(0,T;H^{2 \theta -1 }(\tore)^{3 \times 3})$, it converges weakly, up to a
subsequence, to some $\psi  $ in  $ L^4(0,T;H^{2\theta -1 }(\tore)^{3 \times 3})$. The result above  and uniqueness of the limit,
allows us to claim that $ \psi =\overline{D_{N,\theta}(\bw) \otimes  D_{N,\theta} ({\bw})}$. Consequently
 \begin{align}
     \overline{D_{N,\theta}({\bw}^{n}) \otimes D_{N,\theta}({\bw}^{n})} &\rightharpoonup \overline{D_{N,\theta}(\bw) \otimes D_{N,\theta}({\bw})} &&\textrm{weakly in }
 L^4(0,T;H^{2 \theta - 1 }(\tore)^{3 \times 3}), \label{c22primen}\    
  \end{align}
  We also observe that 
   \begin{align}
     {D_{N,\theta}({\bw}^{n}) \otimes \vec{B}^{n}} &\rightharpoonup {D_{N,\theta}(\bw) \otimes  \vec{B}} &&\textrm{weakly in }
 L^2(0,T;L^{{2}}(\tore)^{3 \times 3}), \label{cb22primen}\\
  {\vec{B}^{n} \otimes D_{N,\theta}({\bw}^{n})} &\rightharpoonup { \vec{B} \otimes D_{N,\theta}(\bw) } &&\textrm{weakly in }
 L^2(0,T;L^{{2}}(\tore)^{3 \times 3}), \label{bc22primenbc22primen}\\
 \overline{\vec{B}^{n} \otimes \vec{B}^{n} } &\rightharpoonup \overline{\vec{B} \otimes \vec{B}} &&\textrm{weakly in }
 L^4(0,T;H^{2 \theta - 1 }(\tore)^{3 \times 3}), \label{bc22primen}
  \end{align}

The above established convergences are clearly sufficient for taking the limit in (\ref{nuw1g})-(\ref{nuweakb1galerkine}) and for concluding that   $ \bw$ , $\vec{B}$ and $ q$   satisfy (\ref{nuweak1})-(\ref{nuweakb1}). 
Moreover, 
from (\ref{nuc122}) and (\ref{nunc322}) and by using that  $ \theta  \ge \frac{5}{6}$, we obtain 
that $ \bw \in  L^2(0,T;\vec{H}^{{\theta}}_{\sigma })$  and $ \partial_t\bw \in  L^2(0,T;\vec{H}^{-\theta}_{ })$ thus 
  we  deduce by a classical argument (\cite{LiMa68, sim87})  that 
 \begin{equation}
 \bw \in  \mathcal{C}(0,T;\vec{L}^{2}_{\sigma }).
\end{equation}
Similarly, we deduce from (\ref{nubnc22}) and (\ref{nubnc322}) that 
\begin{equation}
 \vec{B} \in  \mathcal{C}(0,T;\vec{L}^{{2}}_{\sigma }).
\end{equation}
Furthermore, from  the  strong continuity of $\bw$  and $ \vec{B}$ with respect to the time with value in  $\vec{L}^{{2}}_{\sigma }$,  we deduce   that $\bw(0)=\bw_0$ and $\vec{B}(0)=\vec{B}_{0}$.\\

\textbf{Step 1-iv} (Limit $n \rightarrow \infty$ part II)
We recall that, $ \partial_{t} \vec{w}^n \in L^{\frac{6}{5}}(0,T ; \vec{H}^{2 \theta - \frac{5}{6}}) \hookrightarrow L^{\frac{6}{5}}(0,T ; \vec{H}^{ \theta})$ for any $\theta \ge \frac{5}{6}$, this directly implies that  $ \vec{w}^n \in \mathcal{C}(0,T;\vec{H}^{\theta}_{\sigma }) $ for  almost all $t \in [0,T]$.
It is possible to show that the subsequence $ \left\{\vec{w}^n \right\}_{n \in \N} $ converge strongly to $ \vec{w}$ in   $\mathcal{C}(0,T;\vec{H}^{\theta}_{\sigma }) $ and thus $ \vec{w} \in  \mathcal{C}(0,T;\vec{H}^{\theta}_{\sigma }) $ provided that $\theta \ge  \frac{5}{6}.$  Indeed, it is sufficient to show that  $ \left\{\vec{w}^n\right\}_{n \in \N} $ is a Cauchy sequence in $\mathcal{C}(0,T;\vec{H}^{\theta}_{\sigma })$.  
%  We know that $\vec{w}_{}^{n}  \in L^\infty (0,T;\vec{H}^{\theta})$ and $   {\vec{B}}^{n} \in L^{2}(0,T;\vec{H}^{1}_{\sigma}) $    for $n \in \N$. 
% Let us  define $\delta_{m,n}$ as 
%$$\delta_{m,n}:=|\vec{u}_{\alpha}^{m+n}  - \vec{u}_{}^{m}|_{}^{2}+ \nu \| \vec{u}_{\alpha}^{m+n}- \vec{u}_{}^{m} \|_{}^{2} . $$
The difference $\vec{w}_{}^{n+m}  - \vec{w}_{}^{n} $, for $ m , n \in \N,$  satisfies 
 \begin{equation}
\label{difereence alpha ns}
\begin{array} {llll} \displaystyle
 \frac{d (\vec{w}_{}^{n+m}  - \vec{w}_{}^{n})}{d t}+  \nabla \cdot \overline{(D_{N,\theta}(\vec{w}^{n+m})\otimes(D_{N,\theta}(\vec{w}^{n+m}))} -  \nabla \cdot \overline{(D_{N,\theta}(\vec{w}^{n})\otimes(D_{N,\theta}(\vec{w}^{n})})\\
%\nabla \cdot\left(\overline{(D_{N,\theta}(\vec{w}^{n+m}) - D_{N,\theta}(\vec{w}^{n})) \otimes D_{N,\theta}(\vec{w}^{n})}\right) \\ 
%\hskip 1cm + \nabla \cdot \left(\overline{D_{N,\theta}(\vec{w}^{n+m}) \otimes (D_{N,\theta}(\vec{w}^{n+m}) - D_{N,\theta}(\vec{w}^{n}))}\right)  
 \hskip 2cm =  \nabla \cdot \overline{(\vec{B}^{n+m}\otimes\vec{B}^{n+m})} -  \nabla \cdot \overline{(\vec{B}^n\otimes\vec{B}^n)}.
\end{array}
\end{equation}
%We take $ {\vec{f}^{}}=0$ for simplicity, our problem comes from the nonlinear term.\\
By taking $\mathbb{A}^{}_{\theta} D_{N,\theta}(\vec{w}_{}^{n+m} - \vec{u}_{}^{n})$ as test function in (\ref{difereence alpha ns}) we get 
 \begin{equation}
\label{d alpha ns}
\begin{array} {llll} \displaystyle
 \displaystyle \frac{1}{2}\frac{d}{d t}\|\mathbb{A}^{\frac{1}{2}}_{\theta} D_{N,\theta}^{\frac{1}{2}}(\vec{w}_{}^{n+m}  - \vec{w}_{}^{n})\|_{2}^{2} \le
 \displaystyle \left| \int_{\tore}  {\vec{B}^{n+m}\otimes\vec{B}^{n+m}} -  {\vec{B}^n\otimes\vec{B}^n} : \nabla  D_{N,\theta}(\vec{w}^{n+m}) - D_{N,\theta}(\vec{w}^{n}) \right|\\
  \hskip 2cm
 \displaystyle + \left| \int_{\tore}(D_{N,\theta}(\vec{w}^{n+m}) - D_{N,\theta}(\vec{w}^{n})) \otimes D_{N,\theta}(\vec{w}^{n}) : \nabla (D_{N,\theta}(\vec{w}^{n+m}) - D_{N,\theta}(\vec{w}^{n}))  \right|
 %\\
 %&\le& |(B( (\overline{\vec{u}_{}^{m+n}} - \overline{\vec{u}_{}^{m}}),\vec{u}_{}^{m+n} - \vec{u}_{}^{m} ),\vec{u}_{}^{m}) |\\
 %&\le&| \displaystyle \int_{\mathbb{T}_3}(\overline{\vec{u}_{}^{m+n}} - \overline{\vec{u}_{}^{m}}) \vec{u}_{}^{m} (\nabla \vec{u}_{}^{m+n} - \nabla \vec{u}_{}^{m} )|
 \end{array}
 \end{equation}
 % From Young inequality we have 
 %From 
% Lemma \ref{B:prop 2.3}  combined with Young inequality give
 %we have 
 From the H\"{o}lder inequality combined with Sobolev injection we get, for any $\theta \ge \frac{5}{6}$, that 
  \begin{equation}
\label{di alpha ns}
\begin{array} {llll} 
\displaystyle  \left| \int_{\tore}( D_{N,\theta}({\vec{w}_{}}^{n+m}) - D_{N,\theta}({\vec{w}_{}}^{n})) \otimes D_{N,\theta}(\vec{w}^{n}) : \nabla (D_{N,\theta}(\vec{w}^{n+m}) - D_{N,\theta}(\vec{w}_{}^{n})) \right|\\
 %|\int_{\mathbb{T}_3}(\overline{\vec{u}_{}}^{m+n} - \overline{\vec{u}_{}}^{m}) \vec{u}_{}^{m} (\nabla \vec{u}_{}^{m+n} - \nabla \vec{u}_{}^{m} )| 
 %&\le&\displaystyle \frac{1}{2\nu}  \left| ( \overline{\vec{u}_{}^{m+n}} - \overline{\vec{u}_{}^{m}})\vec{u}_{}^{m} \right|_{}^{2} +\displaystyle \frac{\nu}{2} \left| \nabla \vec{u}_{}^{m+n} - \nabla \vec{u}_{}^{m} \right|_{}^{2} 
% \\
\hskip 3cm \le  \displaystyle C(N) \| \overline{\vec{w}_{}}^{n+m} - \overline{\vec{w}_{}}^{n}\|_{\theta,2 }^{2}   \|\vec{w}_{}^{n} \|_{\theta,2}^{}.
 \end{array}
 \end{equation}
  H\"{o}lder inequality combined with Sobolev injection and the  Young inequality gives
 %we have 
  \begin{equation}
\label{di alpha nsforce}
\begin{array} {llll} 
 \displaystyle \left| \int_{\tore}  {\vec{B}^{n+m}\otimes\vec{B}^{n+m}} -  {\vec{B}^n\otimes\vec{B}^n} : \nabla  D_{N,\theta}(\vec{w}_{}^{n+m}) - D_{N,\theta}(\vec{w}_{}^{n}) \right| \\
 %|\int_{\mathbb{T}_3}(\overline{\vec{u}_{}}^{m+n} - \overline{\vec{u}_{}}^{m}) \vec{u}_{}^{m} (\nabla \vec{u}_{}^{m+n} - \nabla \vec{u}_{}^{m} )| 
 %&\le&\displaystyle \frac{1}{2\nu}  \left| ( \overline{\vec{u}_{}^{m+n}} - \overline{\vec{u}_{}^{m}})\vec{u}_{}^{m} \right|_{}^{2} +\displaystyle \frac{\nu}{2} \left| \nabla \vec{u}_{}^{m+n} - \nabla \vec{u}_{}^{m} \right|_{}^{2} 
% \\
\hskip 2cm   \le {C(N)} \| {\vec{B}^{n+m}\otimes\vec{B}^{n+m}} -  {\vec{B}^n\otimes\vec{B}^n}\|_{1-\theta,2}^{2} 
+ \frac{1}{2} \|    \vec{w}_{}^{n+m} - \vec{w}_{}^{n} \|_{\theta,2}^{2}. 
 \end{array}
 \end{equation}
  Thus we get 
 \begin{equation}
\label{dif alpha ns}
\begin{array} {llll} \displaystyle
 \frac{d}{d t} \|\vec{w}_{}^{n+m} - \vec{w}^{n} \|_{\theta,2}^{2}  \le \displaystyle   C(N,\alpha) \| \vec{w}^{n+m} -\vec{w}_{}^{n}\|_{\theta,2}^{2}  \left( 2 \|\vec{w}_{}^{n}\|_{\theta,2}^{2} + 1 \right)\\ 
 \hskip 5 cm +\displaystyle  C(N, \alpha )\| {\vec{B}^{n+m}\otimes\vec{B}^{n+m}} -  {\vec{B}^n\otimes\vec{B}^n}\|_{1-\theta,2}^{2}
 \end{array}
 \end{equation}
 %One has by Gronwall inequality that 
 By Gr\"{o}nwall inequality we get  
  \begin{equation}
\begin{array} {llll} \displaystyle
 \|\vec{w}_{}^{n+m} - \vec{w}^{n} \|_{\theta,2}^{2} \\
  \hskip 1cm \le \left(\|\vec{w}_{}^{n+m}(0) - \vec{w}_{}^{n}(0) \|_{\theta,2}^{2}  +\displaystyle C(N,\alpha) \int_{0}^{T} \| {\vec{B}^{n+m}\otimes\vec{B}^{n+m}} -  {\vec{B}^n\otimes\vec{B}^n}\|_{1-\theta,2}^{2} dt  \right)\\
 \hskip 2cm \times \exp{ \displaystyle C(N,\alpha) \int_{0}^{T} \left( 2\|\vec{w}_{}^{n} \|_{\theta,2}^{}+ 1 \right) dt} 
   \end{array}
 \end{equation}

 We know that $\vec{w}_{}^{n}  \in L^{\infty} (0,T;\vec{H}^{\theta}_{\sigma })$, thus there exists $C(\alpha, T) \ge 0$ such that $$\exp{\int_{0}^{T} C(N,\alpha) \left( 2 \| \vec{w}_{}^{n}\|_{\theta,2}^{} + 1 \right) dt}  \le C(N,\alpha,T).$$

 We observe that  
 \begin{equation}
\label{dif alpha nsfff}
\begin{array} {llll} \displaystyle
  \int_{0}^{T}\| {\vec{B}^{n+m}\otimes\vec{B}^{n+m}} -  {\vec{B}^n\otimes\vec{B}^n}\|_{1-\theta,2}^{2} dt \le   \displaystyle {\int_{0}^{T} \| {\vec{B}^{n+m}\otimes\vec{B}^{n+m}} -  {\vec{B}\otimes\vec{B}}\|_{1-\theta,2}^{2}}\\
  \hskip 5cm+  \displaystyle \int_{0}^{T} \| \vec{B}^{}\otimes\vec{B}^{} -  \vec{B}^n\otimes\vec{B}^n\|_{1-\theta,2}
  \end{array}
 \end{equation}

  Since for $ n $ tends to  $\infty$, $ \vec{B}^{n}$  converges to $\vec{B}$ strongly in $L^{2}(0,T;\vec{L}^{2}_{\sigma })$ and weakly in $L^{2}(0,T;\vec{H}^{1}_{\sigma })$, we deduce by interpolation that $ \vec{B}^{n}$  converges to $\vec{B}$ strongly in $L^{2}(0,T;\vec{H}^{1-\epsilon }_{\sigma })$, for some $ \epsilon >0$. Thus we derive that 
 $ {\vec{B}^n\otimes\vec{B}^n}$ converges strongly to ${\vec{B}\otimes\vec{B}}$ in $L^{1}(0,T;\vec{H}^{\frac{1}{2}-\epsilon }_{\sigma })$. As $1- \theta \le \frac{1}{6}$ and $\epsilon$ can be chosen  such that $1- \theta < \frac{1}{2}-\epsilon $ we conclude that the right hand side of (\ref{dif alpha nsfff}) tends to zero when $n$ tends to $\infty$.

  Since $\vec{w}_{0} \in \vec{H}^{\theta}_{\sigma}  $ 
 then $\|\vec{w}^{n+m}(0) - \vec{w}^{n}(0) \|_{\theta,2}^{2}$
converges to zero when $n$ goes to $\infty$.
  We deduce  that $\vec{w}_{}^{n+m} - \vec{w}_{}^{n}$ tends to zero in $ \mathcal{C}(0,T;\vec{H}^{\theta}_{\sigma})$. 
 This implies that 
    $\left\{\vec{w}_{}^{n}\right\}_{n \in \N}$ is  a Cauchy sequence in $\mathcal{C}(0,T;\vec{H}^{\theta}_{\sigma })$.\\

\textbf{Step 2} (Existence of regular weak solutions)

Now, we show the existence of regular weak solutions, assuming in addition that $ \vec{B}_0 \in \vec{H}^1_{\sigma}.$

\textbf{Step 2-i} (New uniform estimates via a regularity approach)

Multiplying the $r$th equation in (\ref{nuw1g}) with $\alpha^{2 \theta}|\bk|^{4 \theta }\widehat{D_N}{\bc}^n_{r}(t)+|\bk|^{2 \theta }\widehat{D_N}{\bc}^n_{r}(t)$, and 
multiplying the $r$th equation in (\ref{nuweakb1galerkine}) with $|\vec{k}|^2 {\vec{d}}^n_{r}(t)$, summing the resulted equations over ${r=1},2,...,n$,(in short, we use  $ \mathbb{A}^{}_{\theta} (-\Delta)^{\theta} D_{N,\theta}(\vec{w}^n) $ as  a test function in (\ref{nuw1g}) and  $ - \Delta \vec{B}^{n}$  as  a test function in (\ref{nuweakb1galerkine})), and integrating by parts, we come to the following inequality
 
 \begin{equation}
\begin{array}{lllll}
 \label{apriori1bnu}
\displaystyle \frac{1}{2}  \frac{d}{dt} \left( \| \mathbb{A}^{\frac{1}{2}}_{\theta}D_{N,\theta}^{\frac{1}{2}}(\vec{w}^n)\|_{\theta,2}^2  + \| \nabla {\vec{B}}^n \|_{2}^2 \right)
+ \displaystyle      \mu \| \nabla {\vec{B}}^n \|_{1,2}^2\\
 \hskip 0.5cm  \displaystyle \le \left| \int_{\tore}  (D_{N,\theta}(\vec{w}^n)\cdot \nabla) \vec{B}^{n}    \Delta  \vec{B}^{n}   \right| d\vec{x} +  \displaystyle  \left| \int_{\tore}  (D_{N,\theta}(\vec{w}^n)\cdot \nabla) D_{N,\theta}(\vec{w}^n)  (-\Delta)^{\theta} D_{N,\theta}(\vec{w}^n) \right| d\vec{x}  \\
\hskip 1.5cm    \displaystyle + \left| \int_{\tore}  (\vec{B}^n\cdot \nabla) \vec{B}^{n}  (-\Delta)^{\theta} D_{N,\theta}(\vec{w}^n)     \right| d\vec{x}  + \left| \int_{\tore}  (\vec{B}^{n} \cdot \nabla) D_{N,\theta}(\vec{w}^n)    \Delta  \vec{B}^{n}   \right| d\vec{x},    
 \end{array}
\end{equation}

The first term in right hand side is estimated by 

 \begin{equation}
\begin{array}{lllll}
 \label{apriori1bnuapriori1bnu}
\displaystyle  \left| \int_{\tore}  (D_{N,\theta}(\vec{w}^n)\cdot \nabla) \vec{B}^{n}    \Delta  \vec{B}^{n}   \right| d\vec{x}  \le 
 \|D_{N,\theta}(\vec{w}^n) \|_{\theta,2} \| \nabla\vec{B}^{n}\|_{2}^{\theta - \frac{1}{2}}    \|\Delta  \vec{B}^{n}\|_{2}^{\frac{5}{2}- \theta}
 \end{array}
\end{equation}
Where we have used the Sobolev embedding together with the H\"{o}lder inequality
and the interpolation inequality between $ \vec{L}^2$ and $ \vec{H}^1$.

Similarly, we can estimate the second term in the right hand side  by

 \begin{equation}
\begin{array}{lllll}
 \label{apriori1bnuapriori1bnuapriori1bnu}
\displaystyle    \left| \int_{\tore}  (D_{N,\theta}(\vec{w}^n)\cdot \nabla) D_{N,\theta}(\vec{w}^n)  (-\Delta)^{\theta} D_{N,\theta}(\vec{w}^n) \right| d\vec{x} \le \|  D_{N,\theta}(\vec{w}^n) \|_{\theta,2}   \|  D_{N,\theta}(\vec{w}^n) \|_{\frac{5}{2} - 2\theta,2}  \|  D_{N,\theta} (\vec{w}^n)\|_{2\theta,2}, 
 \end{array}
\end{equation}
since $ \theta \ge \frac{5}{6} $,  we have $\frac{5}{2} - 2\theta \le \theta $  and    $\|  D_{N,\theta}(\vec{w}^n) \|_{\frac{5}{2} - 2\theta} \le \|  D_{N,\theta}(\vec{w}^n) \|_{\theta ,2}$  thus we get 
\begin{equation}
\begin{array}{lllll}
 \label{apriori1bnuapriori1bnuapriori1bnuapriori1bnu}
\displaystyle    \left| \int_{\tore}  (D_{N,\theta}(\vec{w}^n)\cdot \nabla) D_{N,\theta}(\vec{w}^n)  (-\Delta)^{\theta} D_{N,\theta}(\vec{w}^n) \right| d\vec{x} \le \|  D_{N,\theta}(\vec{w}^n) \|_{\theta,2}   \|  D_{N,\theta} (\vec{w}^n)\|_{2\theta,2}^2, 
 \end{array}
\end{equation}

The third term  in the right hand side is estimated by

 \begin{equation}
\begin{array}{lllll}
 \label{apriori1bnuapriori1bnuapriori1bnuapriori1bnuapriori1bnu}
\displaystyle  \displaystyle  \left| \int_{\tore}  (\vec{B}^n\cdot \nabla) \vec{B}^{n}  (-\Delta)^{\theta} D_{N,\theta}(\vec{w}^n)     \right| d\vec{x}    \le \| \vec{B}^n  \|_{1,2}   \|    \vec{B}^{n}  \|_{2,2}   \|  D_{N,\theta}(\vec{w}^n) \|_{2\theta,2}
 \end{array}
\end{equation}

The last term  in the right hand side is estimated by
\begin{equation}
\begin{array}{lllll}
 \label{apriori1bnuapriori1bnuapriori1bnuapriori1bnuapriori1bnuapriori1bnu}
\displaystyle  \left| \int_{\tore}  (\vec{B}^{n} \cdot \nabla) D_{N,\theta}(\vec{w}^n)    \Delta  \vec{B}^{n}   \right| d\vec{x}  \le \| \vec{B}^n\|_{1,2}   \|D_{N,\theta}(\vec{w}^n)\|_{2\theta,2}   \|\Delta \vec{B}^n \|_{2}
 \end{array}
\end{equation}

%  \begin{equation}
% \begin{array}{lllll}
%  \label{apriori1bnu}
% \displaystyle  \left| \int_{\tore}  (\vec{B}^{n} \cdot \nabla) D_{N,\theta}(\vec{w}^n)    \Delta  \vec{B}^{n}   \right| d\vec{x}  \le \| \vec{B}^n\|_{\frac{5}{2}- 2\theta,2}   \|D_{N,\theta}(\vec{w}^n)\|_{2\theta,2}   \|\Delta \vec{B}^n \|_{2}
%  \end{array}
% \end{equation}
Where we have used that   $\|  \vec{B}^n \|_{\frac{5}{2} - 2\theta} \le \|  \vec{B}^n \|_{1 ,2}$ for all $ \theta \ge \frac{5}{6} $.\\

Therefore, 
\begin{equation}
\begin{array}{lllll}
 \label{apriori1bnuapriori1bnuapriori1bnuapriori1bnuapriori1bnuapriori1bnuapriori1bnu111}
\displaystyle \frac{1}{2}  \frac{d}{dt} \left( \| \mathbb{A}^{\frac{1}{2}}_{\theta}D_{N,\theta}^{\frac{1}{2}}(\vec{w}^n)\|_{\theta,2}^2  + \| \nabla {\vec{B}}^n \|_{2}^2 \right) +
 \displaystyle      \mu \| \nabla {\vec{B}}^n \|_{1,2}^2\\
 \hskip 0.5cm     \le  \|D_{N,\theta}(\vec{w}^n) \|_{\theta,2} \| \nabla\vec{B}^{n}\|_{2}^{\theta - \frac{1}{2}}    \|\Delta  \vec{B}^{n}\|_{2}^{\frac{5}{2}- \theta} + \|  D_{N,\theta}(\vec{w}^n) \|_{\theta,2}   \|  D_{N,\theta} (\vec{w}^n)\|_{2\theta,2}^2   \\
\hskip 1.5cm  +  2 \| \vec{B}^n\|_{1,2}   \|D_{N,\theta}(\vec{w}^n)\|_{2\theta,2}   \|\Delta \vec{B}^n \|_{2}     
 \end{array}
\end{equation}
Thus, by using the Young inequality combined with  Lemma \ref{lemme1} we get 

\begin{equation}
\begin{array}{lllll}
 \label{apriori1bnuapriori1bnuapriori1bnuapriori1bnuapriori1bnuapriori1bnuapriori1bnu222}
\displaystyle \ \frac{d}{dt} \left( \alpha \|    \vec{w}^n\|_{2\theta,2}^2  + \|  {\vec{B}}^n \|_{1,2}^2 \right) +
 \displaystyle      \mu \| \nabla {\vec{B}}^n \|_{1,2}^2\\
 \hskip 0.5cm     \le  C(N+1)^5 \|\vec{w}^n \|_{\theta,2}^{\frac{4}{2\theta -1}} \| \nabla\vec{B}^{n}\|_{2}^{2}    + C(N+1)^{3}\|  \vec{w}^n \|_{\theta,2}   \|   \vec{w}^n\|_{2\theta,2}^2   \\
\hskip 1.5cm  +  C (N+1)^2 \| \vec{B}^n\|_{1,2}^2   \|\vec{w}^n\|_{2\theta,2}^2       
 \end{array}
\end{equation}

Therefore
\begin{equation}
\begin{array}{lllll}
 \label{apriori1bnuapriori1bnuapriori1bnuapriori1bnuapriori1bnuapriori1bnuapriori1bnu}
\displaystyle \ \frac{d}{dt} \left( \alpha \|    \vec{w}^n\|_{2\theta,2}^2  + \|  {\vec{B}}^n \|_{1,2}^2 \right) +
 \displaystyle      \mu \| \nabla {\vec{B}}^n \|_{1,2}^2\\
 \hskip 0.5cm   \le C(t) (N+1)^5 \left( \|   \vec{w}^n\|_{2\theta,2}^2   +  \| \vec{B}^{n}\|_{1,2}^{2} \right),         
 \end{array}
\end{equation}

where $C(t)= \left(   \|\vec{w}^n \|_{\theta,2}^{\frac{4}{2\theta -1}} +  \|  \vec{w}^n \|_{\theta,2} +     \| \vec{B}^n\|_{1,2}^2    \right) \in L^1(0,T).$ 

Thus, by using Gronwall's lemma we obtain

\begin{equation}
\label{apriori1bnugron} \sup_{t \in [0,T]}\left( \alpha \|    \vec{w}^n\|_{2\theta,2}^2  + \|  {\vec{B}}^n \|_{1,2}^2 \right) \le   C(\alpha,N)\left( \|    \vec{w}_0\|_{2\theta,2}^2  + \|  {\vec{B}}_0 \|_{1,2}^2 \right)\exp{\int_{0}^{T}C(t)dt}.
\end{equation}

Integrate (\ref{apriori1bnuapriori1bnuapriori1bnuapriori1bnuapriori1bnuapriori1bnuapriori1bnu}) on $[0,T]$ we get 
\begin{equation}
 \label{gronapriori1bnu}
\displaystyle   \sup_{t \in [0,T]}\left( \alpha \|    \vec{w}^n\|_{2\theta,2}^2  + \|  {\vec{B}}^n \|_{1,2}^2 \right)+    \mu \int_{0}^{T}\| \nabla {\vec{B}}^n \|_{1,2}^2 dt
    \le C(\alpha, N, \vec{w}_{0}, \vec{B}_{0},T), 
\end{equation}
where  $C(\alpha, N, \vec{w}_{0}, \vec{B}_{0},T)$ does not depend on $n$.

As a consequence from (\ref{gronapriori1bnu}), we derive that there exists $C>0$  such that 
 \begin{align}
 \label{gronwallnuuniform}
    \| D_N({\bw}^{n}) \|_{L^{\infty}(0,T;\vec{H}^{{2\theta}}_{\sigma})}+ \| {\bw}^{n} \|_{L^{\infty}(0,T;\vec{H}^{{2\theta}}_{\sigma})} +   \| {\vec{B}}^{n} \|_{L^{\infty}(0,T;\vec{H}^{1}_{\sigma})}   \le C \\
    \label{gronwallnuuniform2}
    \| {\vec{B}}^{n} \|_{L^{2}(0,T;\vec{H}^{2}_{\sigma})} \le C
\end{align}

As in the proof of Theorem \ref{TH1}  and by using (\ref{gronwallnuuniform})-(\ref{gronwallnuuniform2}),  we can derive the following estimates
 \begin{align}
 \label{gronwallnuicivbarvbarpressure}
\int_{0}^{T}\|q^{n}\|_{2 \theta + \frac{3}{2} ,2}^2 dt < C, \\
\label{gronwallnuicivtemps}
\int_{0}^{T}  \| \partial_t{\bw^{n}}\|_{2 \theta + \frac{1}{2} ,2}^2  dt < C,\\
\label{gronwallnuicivbtemps}
\int_{0}^{T}  \|\partial_t{\vec{B}^{n}}\|_{2}^2  dt < C. 
\end{align}

\textbf{Step 2-ii} (Limit $n \rightarrow \infty$ part III)

 From {Step 1-iii}, we already know
that there exists a weak solution $( \bw , \vec{B}, q )$   of  (\ref{nuweak1})-(\ref{nuweakb1})  and    not relabeled  subsequence of $( \bw^{n}, \vec{B}^{n},  q^n)$ such that  
\begin{align}
\bw^{n} &\rightarrow \bw &&\textrm{strongly in  }
L^2(0,T;\vec{L}^{{2}}_{\sigma }),
\label{hanuc83icil2}\\
\vec{B}^{n} &\rightarrow \vec{B} &&\textrm{strongly in  }
L^2(0,T;\vec{L}^{{2}}_{\sigma }),
\label{hanubc83icil2}\\
D_{N,\theta}(\bw^{n}) &\rightarrow D_{N,\theta}(\bw) &&\textrm{strongly in  }
L^2(0,T;\vec{L}^{2}_{\sigma }),
\label{hanuDNc83icil2}\\
q^{n}&\rightharpoonup q &&\textrm{weakly in } L^{2}(0,T;H^{ 2\theta -{1} }(\tore)), \label{haweaknuc32}
%\overline{\bw}^n &\rightarrow \overline{\bw} &&\textrm{strongly in  }
%L^q(0,T;L^q(\tore)^3) \textrm{  for all } q< 5 .\label{c82pppp}
\end{align}

 Thus, it follows from the estimates (\ref{gronwallnuuniform})-(\ref{gronwallnuicivbtemps}), the Aubin-Lions compactness lemma
(see \cite{sim87} for example) and the uniqueness of the limit that there are a  not relabeled  subsequence of $(\bw^{n}, \vec{B}^{n}, q^{n})$  and a triplet $(\bw, \vec{B}, q)$ such that
\begin{align}
\bw^{n} &\rightharpoonup^* \bw &&\textrm{weakly$^*$ in } L^{\infty}
(0,T;\vec{H}^{2\theta}_{\sigma }), \label{hanuc122}\\
\vec{B}^{n} &\rightharpoonup^* \vec{B} &&\textrm{weakly$^*$ in } L^{\infty}
(0,T;\vec{H}^{{1}}_{\sigma }), \label{hanubc122}\\
% D_N (\bw^{n}) &\rightharpoonup^* D_N(\bw) &&\textrm{weakly$^*$ in } L^{\infty}
% (0,T;\vec{H}^{2\theta}_{\sigma }), \label{hanuDNc122}\\
% \bw^n &\rightharpoonup \bw &&\textrm{weakly in }
% L^2(0,T;\vec{H}^{\frac{3}{2}}_{\sigma }), \label{nc22}\\
\vec{B}^{n} &\rightharpoonup \vec{B} &&\textrm{weakly in }
L^2(0,T;\vec{H}^{{2}}_{\sigma }), \label{hanubnc22}\\
% D_N\bw^n &\rightharpoonup D_N\bw &&\textrm{weakly in }
% L^2(0,T;\vec{H}^{\frac{3}{2}}_{\sigma }), \label{DNnc22}\\
\partial_t{\bw^{n}}&\rightharpoonup \partial_t\bw &&\textrm{weakly in } L^{2}
(0,T;\vec{H}^{2\theta +\frac{1}{2}}),
\label{hanunc322}\\
\partial_t{\vec{B}^{n}}&\rightharpoonup \partial_t\vec{B} &&\textrm{weakly in } L^{2}
(0,T;\vec{L}^{2}),
\label{haanubnc322}\\
q^{n}&\rightharpoonup q &&\textrm{weakly in } L^{2}(0,T;H^{ 2\theta +\frac{3}{2} }(\tore)), \label{hanuc32}\\
%\bw^n &\rightharpoonup \bw &&\textrm{weakly in } L^{\frac{8}{3}}
%(0,T;L^{\frac{8}{3}}(\partial \Omega)^3). \label{c5.12}\\
\bw^{n} &\rightarrow \bw &&\textrm{strongly in  }
L^2(0,T;\vec{H}^{{\theta}}_{\sigma }),
\label{hanuc83icil2hanuc83icil2}\\
\vec{B}^{n} &\rightarrow \vec{B} &&\textrm{strongly in  }
L^2(0,T;\vec{H}^{{1}}_{\sigma }),
\label{hanubc83icil2hanubc83icil2}\\
D_{N,\theta}(\bw^{n}) &\rightarrow D_{N,\theta}(\bw) &&\textrm{strongly in  }
L^2(0,T;\vec{H}^{\theta}_{\sigma }),
\label{hanuDNc83icil2hanuDNc83icil2}
%\overline{\bw}^n &\rightarrow \overline{\bw} &&\textrm{strongly in  }
%L^q(0,T;L^q(\tore)^3) \textrm{  for all } q< 5 .\label{c82pppp}
\end{align}
Moreover, 
% by using that  $ \theta  \ge \frac{5}{6}$, 
% we obtain 
% that  $ \bw_{,t} \in  L^\infty (0,T;\vec{H}^{\theta}_{\sigma })$, since $ \bw \in  L^\infty (0,T;\vec{H}^{2\theta}_{\sigma })$
% we  deduce by a classical argument (\cite{LiMa68, sim87}) that 
%  \begin{equation}
%  \bw \in  \mathcal{C}(0,T;\vec{H}^{\theta}_{\sigma }).
% \end{equation}
% Similarly, 
we obtain that $  \vec{B} \in 
L^2(0,T;\vec{H}^{{2}}_{\sigma })  $  and $\partial_t\vec{B} \in  L^{2}
(0,T;\vec{L}^{2}_{}) $ thus we deduce by a classical argument (\cite{LiMa68, sim87}) that  
\begin{equation}
 \vec{B} \in  \mathcal{C}(0,T;\vec{H}^{{1}}_{\sigma }).
\end{equation}

% \begin{Rem}
%  As in Step 1-iv,  it  is also possible  to show that  $ \left\{\vec{w}^n\right\}_{n \in N} $ is a Cauchy sequence in $\mathcal{C}(0,T;\vec{H}^{2\theta}_{\sigma })$.   Thus the subsequence $ \left\{\vec{w}^n\right\}_{n \in N} $ converge strongly to $ \vec{w}$ in   $\mathcal{C}(0,T;\vec{H}^{2\theta}_{\sigma }) $ and $ \vec{w} \in  \mathcal{C}(0,T;\vec{H}^{2\theta}_{\sigma })$ provided that $\theta \ge  \frac{5}{6}.$  
% \end{Rem}

% Let us mention also that $D_N({\bw})+ \alpha^{}(-\Delta)^{\theta}D_N({\bw}) \in L^2(0,T;\vec{L}^{2}_{\sigma }),$ and $\vec{B} \in L^2(0,T;\vec{H}^{{2}}_{\sigma })$ hence  $\mathbb{A}_{\theta}D_{N,\theta}{\bw} $  is a possible  test function in the weak formulation (\ref{nuweak1}) and $\vec{B}$ is a possible  test function in the weak formulation (\ref{nuweakb1}). Thus $ \mathbb{A}^{\frac{1}{2}}_{\theta}D_{N,\theta}^{\frac{1}{2}}{\bw} $ and $\vec{B}$ verify for all $t \in [0,T]$  the following equality   
% \begin{equation}
% \begin{array}{lllll}
%  \label{nuapriorileary1}
% \displaystyle \frac{1}{2} \left(\| \mathbb{A}_{\theta}^{\frac{1}{2}}D_{N,\theta}^{\frac{1}{2}}{\bw} \|_{2}^2  + \| {\vec{B}} \|_{2}^2 \right)
% + \displaystyle \mu  \int_{0}^{t}  \left(    \|{\vec{B}} \|_{1,2}^2\right)  \ ds\\
%  \quad  \quad \quad = 
%  \displaystyle \frac{1}{2}\left( \| \mathbb{A}^{\frac{1}{2}}D_N^{\frac{1}{2}}\overline{\bv}^{}_{0} \|_{2}^2   + \| {\vec{B}}_0 \|_{2}^2 \right).
% \end{array}
% \end{equation}
  
\textbf{Step 3} (Unicity of regular weak solutions)
Next, we will  prove the uniqueness of regular weak solutions
among the class of weak solutions.

%Next, we will show the continuous dependence of the  regular weak solutions on the initial data and in particular the uniqueness.\\
Let $ \theta \ge \frac{5}{6}$  and let $({\bw_1,\vec{B}_1},q_1)$ be a weak solution of (\ref{dalpha ns}) on the interval $[0,T]$, with initial values $\bv_1(0) \in \vec{L}^2_{\sigma}$ $,\vec{B}_1(0) \in \vec{L}^2_{\sigma}$ and $({\bw_2,\vec{B}_2},q_2)$ be a regular weak solution of (\ref{dalpha ns}) on the interval $[0,T]$, with initial values $\bv_1(0) \in \vec{L}^2_{\sigma}$ $,\vec{B}_1(0) \in \vec{H}^1_{\sigma}$ 
% be any two solutions of (\ref{dalpha ns}) on the interval $[0,T]$, with initial values $(\bw_1(0),\vec{B}_1(0))$ and $(\bw_2(0),\vec{B}_2(0))$. 
Let us denote by  $\delta \vec{w}_{} =\bw_2-\bw_1$, by $\delta \vec{B} = \vec{B}_2 - \vec{B}_1$ and by $\delta q = q_2 - q_1$.

Then one has
\begin{equation}
\label{matin1matin1}
\begin{split}
\partial_{t} \delta\vec{w}  + \diver(\overline{D_{N,\theta}(\bw_2) \otimes D_{N,\theta}(\bw_2)})- \diver(\overline{D_{N,\theta}(\bw_1) \otimes D_{N,\theta}(\bw_1}))\\
\quad \quad  - \diver(\overline{\vec{B}_2 \otimes \vec{B}_2})+\diver(\overline{\vec{B}_1 \otimes \vec{B}_1})  + \nabla\delta q =0,\\
\partial_{t} \delta\vec{B}  - \mu
\Delta\delta\vec{B} + \diver({D_{N,\theta}(\bw_2) \otimes \vec{B}_2})- \diver({D_{N,\theta}(\bw_1) \otimes \vec{B}_1})\\
\quad \quad  - \diver({\vec{B}_2 \otimes D_{N,\theta}(\bw_2)})+\diver({\vec{B}_1 \otimes D_{N,\theta}(\bw_1)})    = 0,
\end{split}
\end{equation}
and $\delta\vec{w} = 0$, $\delta\vec{B} =0$  at initial time.

Applying $\mathbb{A}_{\theta}^{\frac{1}{2}}  $  to  the first equation of (\ref{matin1matin1})  we obtain

\begin{equation}
\label{matin1ha}
\begin{split}
\mathbb{A}_{\theta}^{\frac{1}{2}}  \partial_{t} \delta\vec{w} +  \mathbb{A}_{\theta}^{\frac{1}{2}} \diver(\overline{D_{N,\theta}(\bw_2) \otimes D_{N,\theta}(\bw_2)})- \mathbb{A}_{\theta}^{\frac{1}{2}} \diver(\overline{D_{N,\theta}(\bw_1) \otimes D_{N,\theta}(\bw_1}))\\
\quad \quad  - \mathbb{A}_{\theta}^{\frac{1}{2}} \diver(\overline{\vec{B}_2 \otimes \vec{B}_2})+ \mathbb{A}_{\theta}^{\frac{1}{2}} \diver(\overline{\vec{B}_1 \otimes \vec{B}_1})  + \mathbb{A}_{\theta}^{\frac{1}{2}} \nabla\delta q =0.\\
\end{split}
\end{equation}

One can take $ \mathbb{A}_{\theta}^{\frac{1}{2}} D_{N,\theta}(\delta\vec{w}) \in \mathcal{C}(0,T;\vec{L}^{{2}}_{\sigma })$
 as test function  in   (\ref{matin1ha}) and $ \delta\vec{B} \in L^2(0,T;\vec{H}^{{1}}_{\sigma })$
 as test function in the second equations of (\ref{matin1matin1}). Let us mention that, $ \partial_{t} \delta\vec{w} \in L^{\frac{6}{5}}(0,T ; \vec{H}^{2 \theta - \frac{5}{6}}) \hookrightarrow L^{\frac{6}{5}}(0,T ; \vec{H}^{ \theta})$  and $\delta\vec{w} \in \mathcal{C}(0,T;\vec{H}^{{\theta}}_{\sigma })$, for any $\theta \ge \frac{5}{6}$,  thus   Lemma 3.4 in \cite{LariostitiMHD} yields to 
 $$
 \left(  \partial_{t} \mathbb{A}_{\theta}^{\frac{1}{2}} \delta\vec{w},  \mathbb{A}_{\theta}^{\frac{1}{2}} D_{N,\theta}(\delta\vec{w}) \right)= \frac{1}{2}\frac{d}{d}\| \mathbb{A}_{\theta}^{\frac{1}{2}}D_{N,\theta}^{\frac{1}{2}}(\delta\vec{w})\|_{2}^2.$$

Since $ \partial_{t} \delta\vec{B} \in L^2(0,T;\vec{H}^{{-1}}_{ })$, by using Lions-Magenes Lemma \cite{LiMa68}  we may justifiably write 
 $$
 \left\langle \partial_{t} \delta\vec{B},  \delta\vec{B} \right\rangle_{\vec{H}^{-1},\vec{H}^{1}}= \frac{1}{2}\frac{d}{d}\| \delta\vec{B}\|_{2}^2.$$

Now, we proceed in the exact way as in the proof of theorem \ref{TH1} in order to obtain the following equality

\begin{equation}
 \label{nusoirab}
\begin{array}{llll}
 \displaystyle \frac{1}{2}\frac{d}{dt} \left(\|\mathbb{A}^{\frac{1}{2}}D_{N,\theta}^{\frac{1}{2}}({\delta \vec{w}})\|_{2}^{2} + \|{\delta \vec{B}}\|_{2}^{2} \right) + \mu \|{\delta \vec{B}}_{}\|_{1,2}^{2} \\
\hskip 1cm= -\displaystyle \int_{\tore} (D_{N,\theta}( \delta \vec{w}) \cdot   \nabla) D_{N,\theta}(\vec{w}_2) D_{N,\theta}(\delta\vec{w})   + \displaystyle \displaystyle \int_{\tore}  D_{N,\theta}(\delta \vec{w}) \otimes \vec{B}_2 : \nabla \delta\vec{B}  \\
\hskip 1.5cm
+ \displaystyle \int_{\tore}   (\delta \vec{B} \cdot \nabla) \vec{B}_2  D_N(\delta\vec{w})  - \displaystyle  \int_{\tore}  \delta \vec{B} \otimes  D_{N,\theta}(\vec{w}_2) : \nabla \delta\vec{B}.
 \end{array}
\end{equation}

Next, we estimate the four integrals in the right hand side of (\ref{nusoirab}). The estimates are obtained by using  H\"{o}lder inequality, Sobolev embedding theorem, the Young inequality  and Lemma \ref{lemme1}.

\begin{equation}
 \label{nunotuniformzero1}
 \begin{split}
   \left| \int_{\tore} (D_{N,\theta}( \delta \vec{w}) \cdot   \nabla) D_{N,\theta}(\vec{w}_2) D_{N,\theta}(\delta\vec{w})  \right|
  %= \displaystyle \left| \left( D_N(\bw_{2}) - D_N(\bw)_{1} \otimes D_N(\bw)_{1}, \nabla D_N(\delta \vec{w}) \right) \right|\\
 % ( {D_N\bw_{2}\times \nabla \times  {D_N\bw}_{2}}-{D_N\bw_{1}\times \nabla \times {D_N\bw}_{1}}, {D_N\delta \vec{w}} )\\
  \le \displaystyle {C(N+1)^6} \|\delta \vec{w}\|_{\theta,2}^{2} \| {\bw}_{2} \|^{}_{2\theta,2} 
% + \displaystyle \frac{C(\alpha)(N+1)^2}{\nu} \|{ \delta \vec{B}}\otimes   {\vec{B}}_{2} \|_{-{\frac{1}{2}},2}^{2} + \frac{\nu \alpha}{4} \|\nabla {\delta \vec{w}} \|_{\frac{1}{2},2}^2 \\
% + \displaystyle \frac{C}{\mu} \|{D_N \delta \vec{w}}\otimes   {\vec{B}}_{1} \|_{2}^{2} + \frac{\mu}{4} \|\nabla {\delta \vec{B}} \|_{2}^2 \\
% + \displaystyle \frac{C}{\mu} \| \delta {\vec{B}}_{} \otimes {D_N  \vec{w}_{1}}   \|_{2}^{2} + \frac{\mu }{4} \|\nabla {\delta \vec{B}} \|_{2}^2. 
\end{split}
\end{equation}
\begin{equation}
 \label{notuniformzero2notuniformzero2}
 \begin{split}
  \displaystyle \left| \int_{\tore}  D_{N,\theta}(\delta \vec{w}) \otimes \vec{B}_2 : \nabla \delta\vec{B}   \right| 
  %= \displaystyle \left| \left( D_N(\bw_{2}) - D_N(\bw)_{1} \otimes D_N(\bw)_{1}, \nabla D_N(\delta \vec{w}) \right) \right|\\
 % ( {D_N\bw_{2}\times \nabla \times  {D_N\bw}_{2}}-{D_N\bw_{1}\times \nabla \times {D_N\bw}_{1}}, {D_N\delta \vec{w}} )\\
% \le \displaystyle \frac{C(\alpha)(N+1)^2}{\nu} \|{D_N \delta \vec{w}}\otimes   {D_N \bw}_{1} \|_{-{\frac{1}{2}},2}^{2} + \frac{\nu \alpha}{4} \|\nabla {\delta \vec{w}} \|_{\frac{1}{2},2}^2 \\
% + \displaystyle \frac{C(\alpha)(N+1)^2}{\nu} \|{ \delta \vec{B}}\otimes   {\vec{B}}_{1} \|_{-{\frac{1}{2}},2}^{2} + \frac{\nu \alpha}{4} \|\nabla {\delta \vec{w}} \|_{\frac{1}{2},2}^2 \\
% \le  \displaystyle \frac{C}{\mu} \| \delta {\vec{B}}_{} \otimes {D_N  \vec{w}_{1}}   \|_{2}^{2} + \frac{\mu }{4} \|\nabla {\delta \vec{B}} \|_{2}^2. 
\le  \displaystyle \frac{C}{\mu} \|{D_{N,\theta}(\delta \vec{w})}\otimes   {\vec{B}}_{2} \|_{2}^{2} + \frac{\mu}{6} \|\nabla {\delta \vec{B}} \|_{2}^2\\
\hskip 1.5cm  \le  \displaystyle \frac{C(N+1)^2}{\mu} \| \delta \vec{w}\|_{\theta,2}^2 \|{\vec{B}}_{2} \|_{1,2}^{2} + \frac{\mu}{6} \|\nabla {\delta \vec{B}} \|_{2}^2, 
\end{split}
\end{equation}

\begin{equation}
 \label{notuniformzero3notuniformzero3}
 \begin{split}
  \displaystyle \left|  \int_{\tore}   (\delta \vec{B} \cdot \nabla) \vec{B}_2  D_N(\delta\vec{w})  \right|  
  %= \displaystyle \left| \left( D_N(\bw_{2}) - D_N(\bw)_{1} \otimes D_N(\bw)_{1}, \nabla D_N(\delta \vec{w}) \right) \right|\\
 % ( {D_N\bw_{2}\times \nabla \times  {D_N\bw}_{2}}-{D_N\bw_{1}\times \nabla \times {D_N\bw}_{1}}, {D_N\delta \vec{w}} )\\
\le \displaystyle \frac{C(N+1)^2}{\alpha } \|    \delta \vec{B} \|_{1,2} \|\vec{B}_2 \|_{2,2}^{}  \| {\delta \vec{w}} \|_{\theta,2} \\
 \hskip 1.5cm \le \displaystyle \frac{C(N+1)^2}{ \mu }  \|\vec{B}_{2} \|^{2}_{2,2} \| {\delta \vec{w}} \|_{\theta,2}^2    + \frac{\mu }{6} \|\delta \vec{B}\|_{1,2}^{2},
% + \displaystyle \frac{C(\alpha)(N+1)^2}{\nu} \|{ \delta \vec{B}}\otimes   {\vec{B}}_{2} \|_{-{\frac{1}{2}},2}^{2} + \frac{\nu \alpha}{4} \|\nabla {\delta \vec{w}} \|_{\frac{1}{2},2}^2 \\
% + \displaystyle \frac{C}{\mu} \|{D_N \delta \vec{w}}\otimes   {\vec{B}}_{1} \|_{2}^{2} + \frac{\mu}{4} \|\nabla {\delta \vec{B}} \|_{2}^2 \\
% + \displaystyle \frac{C}{\mu} \| \delta {\vec{B}}_{} \otimes {D_N  \vec{w}_{1}}   \|_{2}^{2} + \frac{\mu }{4} \|\nabla {\delta \vec{B}} \|_{2}^2. 
\end{split}
\end{equation}

\begin{equation}
 \label{notuniformzero4notuniformzero4}
 \begin{split}
  \displaystyle \left|    \int_{\tore}  \delta \vec{B} \otimes  D_{N,\theta}(\vec{w}_2) : \nabla \delta\vec{B}  \right|
  %= \displaystyle \left| \left( D_N(\bw_{2}) - D_N(\bw)_{1} \otimes D_N(\bw)_{1}, \nabla D_N(\delta \vec{w}) \right) \right|\\
 % ( {D_N\bw_{2}\times \nabla \times  {D_N\bw}_{2}}-{D_N\bw_{1}\times \nabla \times {D_N\bw}_{1}}, {D_N\delta \vec{w}} )\\
% \le \displaystyle \frac{C(\alpha)(N+1)^2}{\nu} \|{D_N \delta \vec{w}}\otimes   {D_N \bw}_{1} \|_{-{\frac{1}{2}},2}^{2} + \frac{\nu \alpha}{4} \|\nabla {\delta \vec{w}} \|_{\frac{1}{2},2}^2 \\
% + \displaystyle \frac{C(\alpha)(N+1)^2}{\nu} \|{ \delta \vec{B}}\otimes   {\vec{B}}_{1} \|_{-{\frac{1}{2}},2}^{2} + \frac{\nu \alpha}{4} \|\nabla {\delta \vec{w}} \|_{\frac{1}{2},2}^2 \\
% \le  \displaystyle \frac{C}{\mu} \| \delta {\vec{B}}_{} \otimes {D_N  \vec{w}_{1}}   \|_{2}^{2} + \frac{\mu }{4} \|\nabla {\delta \vec{B}} \|_{2}^2. 
 \le \displaystyle \frac{C(N+1)^2}{\mu} \|{\delta \vec{B}}\|_{2}^2 \|  D_{N,\theta}(\vec{w}_2) \|_{2\theta,2}^{2} + \frac{\mu}{6} \| {\delta \vec{B}} \|_{1,2}^2
\end{split}
\end{equation}

By using  (\ref{lemme13d}) we have 
\begin{equation}
\label{notuniformnotuniform}
 \begin{split}
  \displaystyle 
   \frac{1}{2}\frac{d}{dt}\left( \alpha^{} \|{\delta \vec{w}}\|_{\theta,2}^{2} + \|{\delta \vec{B}}\|_{2}^{2} \right)  + \mu \|  {\delta \vec{B}}_{}\|_{1,2}^{2}  \hskip 3cm \\
\le 
  \frac{1}{2}\frac{d}{dt} \left( \|\mathbb{A}^{\frac{1}{2}}_{\theta}D_{N,\theta}^{\frac{1}{2}}({\delta \vec{w}})\|_{2}^{2}+ \|{\delta \vec{B}}\|_{2}^{2}\right)   + \mu \|  {\delta \vec{B}}_{}\|_{1,2}^{2} \end{split}
\end{equation}
From (\ref{soirab})-(\ref{notuniformnotuniform})  we get 
\begin{equation}
 \begin{split}
  \displaystyle 
  \frac{d}{dt}\left(\alpha^{}  \|{\delta \vec{w}}\|_{\theta,2}^{2} + \|{\delta \vec{B}}\|_{2}^{2} \right)  + \mu \|{\delta \vec{B}}_{}\|_{1	,2}^{2}   \hskip 5cm\\
\le \displaystyle {C(N+1)^6} \|\delta \vec{w}\|_{\theta,2}^{2} \| {\bw}_{2} \|^{}_{2\theta,2} +
  \displaystyle \frac{C(N+1)^2}{\mu} \| \delta \vec{w}\|_{\theta,2}^2 \|{\vec{B}}_{2} \|_{1,2}^{2} \\
 \hskip 2cm +\frac{C(N+1)^2}{ \mu}   \| {\delta \vec{w}} \|_{\theta,2}^2  \|\vec{B}_{2} \|^{2}_{2,2}    
+  \displaystyle \frac{C(N+1)^2}{\mu} \|{\delta \vec{B}}\|_{2}^2 \|  \vec{w}_2 \|_{2\theta,2}^{2}.
\end{split}
\end{equation}

Hence, 

\begin{equation}
 \begin{split}
  \displaystyle 
   \frac{d}{dt}\left(\alpha^{}  \|{\delta \vec{w}}\|_{\theta,2}^{2} + \|{\delta \vec{B}}\|_{2}^{2} \right)   + \mu \|{\delta \vec{B}}_{}\|_{1	,2}^{2}  \hskip 6cm \\
\le \displaystyle \frac{C(N,\alpha)}{ \mu} \left( \alpha \|\delta \vec{w}\|_{\theta,2}^{2} + \|\delta \vec{B}\|_{2}^{2} \right) \left( \| {\bw}_{2} \|^{}_{2\theta,2}+ \| {\bw}_{2} \|^{2}_{2\theta,2} + 
  \displaystyle   \|{\vec{B}}_{2} \|_{2,2}^{2} \right).
\end{split}
\end{equation}

 % \begin{equation}
% \begin{split}
%  \displaystyle 
%\|{\delta \vec{w}}_{,t}\|_{2}^{2} +\nu \|\nabla{\delta \vec{w}}_{}\|_{2}^{2} & \le \displaystyle
%\frac{4}{\nu} \|\widetilde{\delta \vec{w}}
%\overline{\bw}_{1} \|_{2}^{2} +  \frac{4}{\nu} \|\widetilde{\vec{v}_2}
%\overline{\delta \vec{w}} \|_{2}^{2}    \\
%&\le \displaystyle \frac{4}{\nu} \|\widetilde{\delta \vec{w}}\|_{\frac{1}{2}-2\theta_2,2}^{2}
%\|\overline{\bw}_{1} \|^{2}_{1+2\theta_2}+  \displaystyle\frac{4}{\nu} \|\overline{\delta \vec{w}}\|_{\frac{1}{2}-2\theta_1,2}^{2}
%\|\widetilde{\bw}_{2} \|^{2}_{1+2\theta_1}  \\ &\le \displaystyle \frac{1}{\alpha^{{2\theta_1}+{2\theta_2}}}\frac{4}{\nu} \|{\delta \vec{w}}\|_{2}^{2}
%\left(\|{\bw}_{1} \|^{2}_{1,2} + 
%\|{\bw}_{2} \|^{2}_{1,2}\right)
%\end{split}
%\end{equation}
%\begin{equation}
% \begin{split}
%  \displaystyle 
%\|{\delta \vec{w}}_{,t}\|_{2}^{2} +\nu \|\nabla{\delta \vec{w}}_{}\|_{2}^{2} & \le \displaystyle
%\frac{4}{\nu} \|\widetilde{\delta \vec{w}}
%\overline{\bw}_{1} \|_{2}^{2} +  \frac{4}{\nu} \|\widetilde{\vec{v}_2}
%\overline{\delta \vec{w}} \|_{2}^{2}    \\
%&\le \displaystyle \frac{4}{\nu} \|\widetilde{\delta \vec{w}}\|_{\frac{3}{4}+2\theta_1-\theta_2,2}^{2}
%\|\overline{\bw}_{1} \|^{2}_{\frac{3}{4}-2\theta_1+\theta_2}+  \displaystyle\frac{4}{\nu} \|\overline{\delta \vec{w}}\|_{\frac{1}{2}-2\theta_1+\theta_2,2}^{2}
%\|\widetilde{\bw}_{2} \|^{2}_{1+2\theta_1-\theta_2}  \\ &\le \displaystyle \frac{1}{\alpha^{{2\theta_1}+{2\theta_2}}}\frac{4}{\nu} \|{\delta \vec{w}}\|_{2}^{2}
%\left(\|{\bw}_{1} \|^{2}_{1,2} + 
%\|{\bw}_{2} \|^{2}_{1,2}\right)
%\end{split}
%\end{equation}
 Since   $\| {\bw}_{2} \|^{}_{2\theta,2} + \| {\bw}_{2} \|^{2}_{2\theta,2} +\|{\vec{B}}_{2} \|_{1,2}^{2} \in L^1([0,T])$, we conclude by using Gronwall's inequality 
the continuous dependence of the solutions on the initial data. In particular, if ${\delta \vec{w}}^{}_{0}= {\delta \vec{B}}^{}_{0}=0$ then ${\delta \vec{w}}= {\delta \vec{B}}=0$ and the solutions are unique for all $t \in [0,T] .$  Since $T>0$ is arbitrary this solution may be uniquely extended for all time.\\ 
This finishes the proof of Theorem \ref{Theorem2}.\\

\section{Limit when $N \rightarrow \infty $ in the double viscous case }
\hskip 0.5cm 
Let $(\vec{w}_N, \vec{B}_N, q_N ) $   be the unique solution of (\ref{dalpha ns}) constructed in theorem (\ref{TH1}) with $ N>0$ fixed and $  \frac{1}{2} \le \theta \le 1$.   The main result of this section is the following.

% In this section, we  take the limit $N \rightarrow \infty $ and look for the
%  optimal value of $\theta$ 
%     in order to show that  $(\vec{w}_N, q_N ) $  converges,  up to subsequences, to a solution of the mean 
% rotational  Navier-Stokes equations. 

% We note that this result holds true  for all $ \theta >0$, however we will give the proof for all $  \theta \in [ \frac{1}{6}, 1[$. 
 \begin{Theorem}
\label{1deuxieme}
Let $ \alpha >0 $ and $   \frac{1}{2} \le \theta < 1$, then from the sequence $    \left\{(\vec{w}_N, \vec{B}_N,  q_N )\right\}_{N \in \N}$, one can extract  a not relabeled subsequence $    \left\{(\vec{w}_N, \vec{B}_N ,q_N )\right\}_{N \in \N} $ such that when $ N$ tends to $  \infty$:\\ 
$(\bw_{N},\vec{B}_N , q_{N}) \rightarrow (\bw, \vec{B}, q)$ 
where 
$$ (\bw, \vec{B}, q) \in L^{\infty}([0,T];\vec{H}^{\theta}_{\sigma})\cap L^2([0,T];\vec{H}^{1+\theta}_{\sigma})
\times L^{\infty}([0,T];\vec{L}^{2}_{\sigma})\cap L^2([0,T];\vec{H}^{1}_{\sigma}) \times L^2([0,T];H^{-\frac{1}{2}+2\theta}(\tore)) $$
is a distributional solution of  the following system with periodic boundary conditions  
\begin{equation}
\label{athetadalpha ns}
 \left\{
 \begin{array} {llll} \displaystyle
  \partial_t \vec{w}^{} + \diver( \overline{ \mathbb{A}_{\theta}(\vec{w})\otimes   \mathbb{A}_{\theta}
(\vec{w}^{}) }) -  \diver( \overline{ \vec{B} \otimes  
\vec{B}^{} }) - \nu \Delta \vec{w}^{} + \nabla
q^{} = 0,\\
\displaystyle \partial_t \vec{B}^{}+ \diver( { \vec{B}\otimes  \mathbb{A}_{\theta}
(\vec{w}^{}) }) -  \diver( \mathbb{A}_{\theta}(\vec{w})\otimes 
\vec{B}^{}) - \mu \Delta \vec{B}^{} = 0,\\
\displaystyle \diver\vec{w}^{} = \diver\vec{B}^{}=0,\
 \displaystyle \int_{\mathbb{T}_3} \vec{w}^{}=  \displaystyle \int_{\mathbb{T}_3} \vec{B}^{}= 0,\\
\displaystyle \vec{w}^{}_{t=0}=\vec{w}_{0}^{}=\overline{\bv_0}, \vec{B}^{}_{t=0}=\vec{B}_{0}^{}.
\end{array}
\right.
\end{equation}

% \begin{align}
% \nabla \cdot  \bw_{} &=0, \label{AthetallBM}\\
% \bw_{,t} -  \overline{\mathbb{A}_{\theta} \bw \times \nabla \times \mathbb{A}_{\theta}\bw} -  \nu \Delta \bw + \nabla q
%  &=  \overline{\bef},\\
% \int_{\tore} \bw  =0, \int_{\tore} q &= 0, \\
% \bw(\bx, 0) = \bw_{0}(\bx)& =\overline{\bv_{0}}. \label{AthetallBLM}
% \end{align}
 The sequence  $ \left\{\bw_{N}\right\}_{N \in \N}$  converges strongly
to $\bw$ in the space $L^2([0,T];\vec{H}^{s})  $ for all  $  s< 1+\theta, $   while the sequence
$ \left\{\vec{B}_{N}\right\}_{N \in \N}$  converges strongly
to $\vec{B}$ in the space $L^2([0,T];\vec{H}^{s})  $ for all  $  s< 1, $ 	and
 $  \left\{q_{N}\right\}_{N \in \N} $  converges weakly to $q$ in the space $ L^2([0,T];H^{-\frac{1}{2}+2\theta}(\tore)).$
 \end{Theorem}

\begin{Rem}
  The restriction in the case of the Approximate Deconvolution Model (ADM)  in \cite{bresslilewandowski} is $\theta > \frac{3}{4}$.
 If we consider the case $ \vec{B} =0$, then equations (\ref{dalpha ns}) with $ \nu >0$
  reduce to   the (ADM). Thus,  the same approach used here can be used to the ADM  to get similar results with $\theta \ge \frac{1}{2}.$
\end{Rem}

% In order to prove Theorem \ref{1deuxieme} we need to
% reconstruct a uniform estimates for $\bw_{N}$. Since some of the a priori  estimates used  in the proof of Theorem \ref{TH1}    dependent   on $N$  we can not use them. Instead we  show in the proof of Theorem \ref{1deuxieme}  that $D_{N,\theta}\bw_{N}$ belongs to $ L^{\infty}([0,T];\vec{L}^{2})\cap L^2([0,T];\vec{H}^{1}_{\sigma })$ uniformly with respect to $N$.\\

Before proving Theorem \ref{1deuxieme}, we first record the following  Lemma.
%On commence par le Lemme suivant:
%We start by the following Lemma: \ref{1deuxieme} 
% On prend $2\phi\vec{u}_{\alpha}$ comme une fonction test dans (\ref{alpha ns}). On note que la condition $\theta \ge 1/4 $ assure que  toutes les termes obtenues sont bien définies, en particulier l'intégrale $$\displaystyle 2\int_{0}^{T} \int_{\tore} \overline{\vec{u}_{\alpha}}\nabla \vec{u}_{\alpha} \cdot \vec{u}_{\alpha} \phi   \  d\vec{x} dt $$ est finie en utilisant le fait que $ \overline{\vec{u}_{\alpha}}\nabla \vec{u}_{\alpha} \in L^2([0,T];\vec{H}^{-1})$ et $2\phi\vec{u}_{\alpha} \in L^2([0,T];\vec{H}^{1})$.\\ 
%Maintenant, des intégrations par parties en utilisant le fait que $\phi(T,\cdot)=\phi(0,\cdot)=0$ et l'identité suivante : 
%\begin{equation}
%2\int_{\tore} \overline{\vec{u}_{\alpha}}\nabla \vec{u}_{\alpha} \cdot \vec{u}_{\alpha} \phi   \  d\vec{x} = \int_{\tore} \overline{\vec{u}_{\alpha}} |\vec{u}_{\alpha}|^2   \cdot \nabla \phi   \  d\vec{x}
%\end{equation}
%montrent que  $(\vec{u}_{\alpha}, p_{\alpha})$ vérifie bien (\ref{local alpha}).\\
%
%Pour prendre la limite $\alpha \rightarrow 0$ on a besoin de montrer que pour tout $\bu_\alpha \in L^p(0,T;L^p(\tore)^3) $, $2 \le p<10/3.$  on a :  
%\begin{align}
%\overline{\bu_\alpha}  &\rightarrow \bu &&\textrm{strongly in  }
%L^p(0,T;L^p(\tore)^3) \textrm{ pour tout }  2 \le p<10/3.
%\end{align}
%Cela est le but de ces deux Lemmes suivants. 
\begin{Lem}
\label{fourierdiscret}
Let $ \theta >0 $ and assume that $\bv  \in L^{2}([0,T], \vec{H}^{2 \theta} )$.   Then   
\begin{align} D_{N,\theta}(\bv)  
  &\rightarrow  \mathbb{A}_{\theta}(\bv)  \textrm{ strongly in }
L^2(0,T;\vec{L}^2),  \textrm{ when } N  \rightarrow \infty, 
 \end{align}
and  there exist a constant C  independent from $N$  such that \begin{align} \| D_{N,\theta}(\bv) \|_{2}  \le \| \mathbb{A}_{\theta} (\bv)\|_{2 }  \le C(\alpha^{\theta}) \|  \bv\|_{2\theta,2 }.    \end{align}
\end{Lem}
\textbf{Proof.} The first part of this lemma is given  in \cite{bresslilewandowski}, see also in \cite{HaniAli}  for  an alternative proof. The second part  is a direct consequence from the  property  (\ref{haniali}) of the operator $D_{N,\theta}$,  the relation (\ref{lemmeatheta2})  and Poincaré inequality.

\textbf{Proof of Theorem \ref{1deuxieme}.}
  The proof of Theorem \ref{1deuxieme} follows the lines of the proof of  Theorem 4.1 in \cite{HaniAli}.
  %that we have to modify in order to treat the cases when $\theta < 1$. 
% The only difference is that instead of  getting a bound on
% the sequence $\partial_t D_{N,\theta} \bw_N $ in a suitable space we use Lemma \ref{fourierdiscret}. 
First, we need to
 reconstruct a uniform estimates for $(\bw_{N}, \vec{B}_N, q_N)$ with respect to $N$.

%De la relation entre $\bu_\alpha$ et $\overline{\bu_\alpha}$ on peut montrer que pour tout $\bu_\alpha \in L^p(0,T;L^p(\tore)^3) $ lorsque $\alpha$ tend vers 0 on a :  
%\begin{align}
%\overline{\bu_\alpha}  &\rightarrow \bu &&\textrm{strongly in  }
%L^p(0,T;L^p(\tore)^3) \textrm{ pour tout } p>1.
%\end{align}
%(Voir Lemme 3.1 page 8)

\textbf{Step 1} (Uniform estimates with respect to N)
%Hereafter we  restrict ourselves to the critical exponent $\theta=\frac{1}{2}$. \\
Following  the proof of Theorem \ref{TH1}  we obtain  
that the solution of (\ref{dalpha ns}) satisfies 
\begin{equation}
\begin{array}{lll}
 \label{undemiapriori12leray}
\displaystyle \frac{1}{2}\left(\| \mathbb{A}_{\theta}^{\frac{1}{2}}D_{N,\theta}^{\frac{1}{2}}({\bw_N}) \|_{2}^2 +\| {\vec{B}_N} \|_{2}^2\right)
+ \displaystyle \int_{0}^{t} \left( \nu \| \mathbb{A}_{\theta}^{\frac{1}{2}}D_{N,\theta}^{\frac{1}{2}}({\bw_N}) \|_{1,2}^2 + \mu\|  {\vec{B}_N} \|_{1,2}^2 \right) \ ds\\
 \quad \le 
 \displaystyle \frac{1}{2} \left(\| \mathbb{A}_{\theta}^{\frac{1}{2}}D_{N,\theta}^{\frac{1}{2}}(\overline{\bv}_{0}) \|_{2}^2 + \| {\vec{B}}_{0} \|_{2}^2 \right),
\end{array}
\end{equation}

 consequently  as $\| \mathbb{A}_{\theta}^{\frac{1}{2}}D_{N,\theta}^{\frac{1}{2}}(\overline{\bv}_{0}) \|_{2}^2 \le  \| {\bv}_{0} \|_{2}^2$ we can bound the right hand side by a constant C which  is independent from $N$.

We deduce from (\ref{undemiapriori12leray}) and Lemma \ref{lemme1} that 
\begin{equation}
\label{vbar1}
 {D_{N,\theta}(\bw_N)} \in L^{\infty}(0,T ; \vec{L}^{2}_{\sigma }) \cap L^{2}(0,T ; \vec{H}^{1}_{\sigma }), \textrm{ uniformly with respect to }  N,
 \end{equation}
and 
\begin{equation}
\label{vbar2}
\begin{split}
{\bw_N} \in L^{\infty}(0,T ; \vec{H}^{\theta}_{\sigma }) \cap L^{2}(0,T ; \vec{H}^{1+\theta}_{\sigma }), \textrm{ uniformly with respect to } N,\\
 {\vec{B}_N} \in L^{\infty}(0,T ; \vec{L}^{2}_{\sigma }) \cap L^{2}(0,T ; \vec{H}^{1}_{\sigma }), \textrm{ uniformly with respect to } N.
 \end{split}
 \end{equation}
 
 We observe from (\ref{vbar1}) that
 \begin{equation}
\label{DNvbarsansbar}
\begin{split}
{D_{N,\theta} ({\bw_N}) \otimes   D_{N,\theta} ({\bw_N})} \in  L^{2}(0,T ;{H}^{-\frac{1}{2}}(\tore)^{3 \times 3 } ),\\
D_{N,\theta} ({\bw_N}) \otimes  \vec{B}_N \in  L^{2}(0,T ;{H}^{-\frac{1}{2}}(\tore)^{3 \times 3 } ),\\
\vec{B}_N \otimes   D_{N,\theta} ({\bw_N}) \in  L^{2}(0,T ;{H}^{-\frac{1}{2}}(\tore)^{3 \times 3 } ),\\
 \vec{B}_N \otimes  \vec{B}_N  \in  L^{2}(0,T ;{H}^{-\frac{1}{2}}(\tore)^{3 \times 3 } ).
\end{split}
\end{equation}
Thus, from (\ref{DNvbarsansbar}) and  (\ref{lemme2.2}) we obtain  
 \begin{equation}
\label{vbarvbar}
\begin{split}
\overline{{D_{N,\theta} ({\bw_N}) \otimes   D_{N,\theta} ({\bw_N})}} \in  L^{2}(0,T ;{H}^{-\frac{1}{2} + 2\theta }(\tore)^{3} ),\\
\overline{{\vec{B}_N} \otimes {\vec{B}_N}} \in  L^{2}(0,T ;{H}^{-\frac{1}{2} + 2\theta }(\tore)^{3 \times 3} ).
\end{split}
\end{equation}
For the pressure term $q_N$, we deduce  that it verifies the following equation 
%Consequently from the Calderon-Zygmund theory the following equation for the pressure 
\begin{equation}
\begin{split}
\displaystyle \Delta^{} q_N=  -\diver \diver  \left(\overline{D_{N,\theta}({\bw}_N) \otimes D_{N,\theta}({\bw}_N)} - \overline{\vec{B}_N \otimes  \vec{B}_N}  \right),
\end{split}\label{pressure}
\end{equation}
consequently, the classical elliptic theory combined with (\ref{vbarvbar}) implies that 
 \begin{equation}
\label{vbarvbarpressure}
\int_{0}^{T}\|q_{N}\|_{-\frac{1}{2} + 2\theta,2}^2 dt < K, \textrm { uniformly with respect to  } N.
\end{equation}
From (\ref{dalpha ns}), (\ref{DNvbarsansbar}) and  (\ref{vbarvbar}) we also obtain that 
 \begin{equation}
\label{vtemps}
\begin{split}
\int_{0}^{T}  \|\partial_t\bw_{N}\|_{-\frac{3}{2} + 2\theta,2}^2  dt < K, \textrm{ uniformly with respect to } N, \\
\int_{0}^{T}  \|\partial_t \vec{B}_{N}\|_{-\frac{3}{2} ,2}^2  dt < K, \textrm{ uniformly with respect to } N.
\end{split}
\end{equation}

\textbf{Step 2} (Passing to the limit $ N \rightarrow \infty$) The central issues is how to take the limit in the nonlinear terms $ {D_{N,\theta} ({\bw_N}) \otimes   D_{N,\theta}({\bw_N})}$,  $ \vec{B}_N \otimes  D_{N,\theta}({\bw_N})$  and $ D_{N,\theta}({\bw_N}) \otimes   \vec{B}_N$.

From the  Aubin-Lions compactness Lemma  (the same arguments as in section 3) we can find a not relabeled  subsequence $\left\{(\vec{w}_{N}, D_{N, \theta}(\vec{w}_{N}), \vec{B}_N, q_{N})\right\}_{N \in \N}$  and $(\vec{w}_{},  \vec{z}, \vec{B}, q_{})$ such that when   $N \rightarrow \infty $  we have: 
%La solution $(\vec{u}_{\alpha_j}, p_{\alpha_j})$  appartient a l'espace d'énergie de Leray donc d'aprés le Lemme de compacité de Aubin-Lions on peut extraire une sous suite tel que lorque $\alpha_j$ tend vers 0 on a: 

\begin{align}
\bw_N &\rightharpoonup^* \bw &&\textrm{weakly$^*$ in } L^{\infty}
(0,T;\vec{H}^{\theta}_{\sigma }), \label{limitNc122}\\
\bw_N &\rightharpoonup \bw &&\textrm{weakly in }
L^2(0,T;\vec{H}^{1+\theta}_{\sigma }), \label{limitNc22}\\
D_{N,\theta}(\bw_N) &\rightharpoonup \vec{z} &&\textrm{weakly in }
L^2(0,T;\vec{H}^{1}_{\sigma }), \label{limitNDNc22}\\
\vec{B}_N &\rightharpoonup^* \vec{B} &&\textrm{weakly$^*$ in } L^{\infty}
(0,T;\vec{L}^{2}_{\sigma }), \label{BlimitNc122}\\
\vec{B}_N &\rightharpoonup \vec{B} &&\textrm{weakly in }
L^2(0,T;\vec{H}^{1}_{\sigma }), \label{BlimitNc22}\\
\partial_t \bw_{N}&\rightharpoonup \partial_t \bw &&\textrm{weakly in } L^{2}
(0,T;\vec{H}^{-\frac{3}{2} + 2\theta}_{}),
\label{limitNc322}\\
\partial_t \vec{B}_{N}&\rightharpoonup \partial_t \vec{B} &&\textrm{weakly in } L^{2}
(0,T;\vec{H}^{-\frac{3}{2}}_{}),
\label{BlimitNc322}\\
q_N&\rightharpoonup q &&\textrm{weakly in } L^{2}(0,T;H^{-\frac{1}{2}+2\theta}
(\tore)), \label{limitNc32}\\
%\bw^n &\rightharpoonup \bw &&\textrm{weakly in } L^{\frac{8}{3}}
%(0,T;L^{\frac{8}{3}}(\partial \Omega)^3). \label{c5.12}\\
\bw_N &\rightarrow \bw &&\textrm{strongly in  }
L^2(0,T;\vec{H}^{s}_{\sigma }) \textrm{ for all } s < 1+\theta,
\label{limitNc83icil2}\\
\vec{B}_N &\rightarrow \vec{B} &&\textrm{strongly in  }
L^2(0,T;\vec{H}^{s}_{\sigma }) \textrm{ for all } s < 1.
\label{BlimitNc83icil2}
%\overline{\bw}^n &\rightarrow \overline{\bw} &&\textrm{strongly in  }
%L^q(0,T;L^q(\tore)^3) \textrm{  for all } q< 5 .\label{c82pppp}
\end{align}

 The goal is
to prove 
 \begin{align}
 D_{N,\theta}({\bw}_N) \otimes \ D_{N,\theta}({\bw}_N) &\rightarrow  \mathbb{A}_{\theta}(\bw) \otimes \mathbb{A}_{\theta}(\bw) &&\textrm{strongly in  }
L^1(0,T;L^{1}(\tore)^{3 \times 3}),\label{DNc82int}\\
 D_{N,\theta}({\bw}_N)  \otimes  {\vec{B}}_N &\rightarrow  \mathbb{A}_{\theta}(\bw) \otimes \vec{B} &&\textrm{strongly in  }
L^1(0,T;L^{1}(\tore)^{3 \times 3}),\label{B1DNc82int}\\
{\vec{B}}_N   \otimes    D_{N,\theta}({\bw}_N) &\rightarrow  \vec{B} \otimes  \mathbb{A}_{\theta}(\bw)  &&\textrm{strongly in  }
L^1(0,T;L^{1}(\tore)^{3 \times 3}).\label{B2DNc82int}
 \end{align}
Thus, it remains to show that 
\begin{align}
 D_{N,\theta}({\bw}_N) &\rightarrow  \mathbb{A}_{\theta} ({\bw})  &&\textrm{strongly in  }
L^2(0,T;\vec{L}^{2}).\label{hhdac82int}
 \end{align}
 In order to show (\ref{hhdac82int}), we compute directly the difference between  $ D_{N,\theta}({\bw}_N) $ and  $\mathbb{A}_{\theta} ({\bw})$ as follows 
\begin{equation}  
 \begin{split}
\| D_{N,\theta}({\bw}_N) - \mathbb{A}_{\theta}({\bw})\|_{L^2(0,T;\vec{L}^{2})} \le \| D_{N,\theta}({\bw}_N) - D_{N,\theta}({\bw})\|_{L^2(0,T;\vec{L}^{2})} + \| D_{N,\theta}({\bw})-  \mathbb{A}_{\theta}({\bw})  \|_{L^2(0,T;\vec{L}^{2})}\\
\le \| \mathbb{A}_{\theta} \left( {\bw}_N - {\bw} \right) \|_{L^2(0,T;\vec{L}^{2})} + \| D_{N,\theta}({\bw})- \mathbb{A}_{\theta}({\bw})\|_{L^2(0,T;\vec{L}^{2})}\\
 \le C\| {\bw}_N -{\bw}\|_{L^2(0,T;\vec{H}^{2\theta})} + \| D_{N,\theta}({\bw})- \mathbb{A}_{\theta}({\bw})\|_{L^2(0,T;\vec{L}^{2})}, 
 \end{split}
 \end{equation}
 where we have used the second part of lemma \ref{fourierdiscret}. 
 Hence, using that $ \theta <1$, (\ref{limitNc22})  and the first part of lemma  (\ref{fourierdiscret}) we deduce (\ref{hhdac82int}). 
 Finally the result above, (\ref{hhdac82int}), combined with (\ref{limitNDNc22})  and the  uniqueness of the limit, allows us to deduce  

 \begin{align}
     D_{N,\theta}({\bw}_N) &\rightharpoonup    \mathbb{A}_{\theta}({\bw}) &&\textrm{weakly in }
 L^2(0,T;\vec{H}^{1}_{\sigma}). \label{c22prime0}   
  \end{align}

Consequently we get (\ref{DNc82int}) and  by using  (\ref{DNvbarsansbar}) we obtain 
 %from 
%(\ref{DNvbarsansbar}), \ref{lemme2.2}   and the  uniqueness of the limit, allows us to deduce 

\begin{align}
     \overline{D_{N,\theta}({\bw}_N) \otimes  D_{N,\theta}({\bw}_N)} &\rightharpoonup \overline{\mathbb{A}_{\theta}(\bw) \otimes   \mathbb{A}_{\theta}({\bw})} &&\textrm{weakly in }
 L^2(0,T;H^{-\frac{1}{2} + 2\theta}(\tore)^{3 \times 3}). \label{c22prime}  
  \end{align}
  
 Similarly we get 
 \begin{align}
     { D_{N,\theta}({\bw}_N) \otimes \vec{B}_N } &\rightharpoonup \mathbb{A}_{\theta}(\bw) \otimes  \vec{B} &&\textrm{weakly in }
 L^2(0,T;H^{-\frac{1}{2}}(\tore)^{3\times 3}), \label{B1c22prime} \\
  {   \vec{B}_N \otimes  D_{N,\theta}({\bw}_N)  } &\rightharpoonup \vec{B}  \otimes  \mathbb{A}_{\theta}(\bw)  &&\textrm{weakly in }
 L^2(0,T;H^{-\frac{1}{2}}(\tore)^{3\times 3}). \label{B2c22prime} 
  \end{align}

These convergence results allow us  to prove   that $ (\bw, \vec{B}, q)$ is a distributional  solution to the mean  MHD equations (\ref{athetadalpha ns}). This finishes the proof of Theorem \ref{1deuxieme}.

\end{document}